\documentclass[10pt]{article}

\usepackage[utf8]{inputenc}
\usepackage{amsmath} 
\usepackage{amssymb} 
\usepackage{amsthm}
\usepackage{xcolor}
\usepackage{parskip}
\usepackage[T1]{fontenc}
\usepackage[french]{babel}

\newtheorem{nota}{Notation}

\newtheorem{proposition}{Proposition}[subsection]
\newtheorem{Remarque}[proposition]{Remarque}
\newtheorem{conjecture}[proposition]{Conjecture}
\newtheorem{fait}[proposition]{Fait}
\newtheorem{definition}[proposition]{Définition}
\newtheorem{theoreme}[proposition]{Théorème}
\newtheorem{Lemme}[proposition]{Lemme}
\newtheorem{corollaire}[proposition]{Corollaire}

\title{Quasi groupes de Frobenius dimensionnels}
\author{Samuel Zamour \footnote{Institut Camille Jordan (ICJ), Université Claude Bernard Lyon 1, e-mail : zamour@math.univ-lyon1.fr}}
\date{11 octobre 2022}

\begin{document}
\maketitle
\begin{abstract}
    Nous nous intéressons à une classe de groupes, les quasi groupes de Frobenius (avec involutions), dont la structure interne généralise celle des groupes classiques $\mbox{GA}_1(\mathbb{C})$, $\mbox{PGL}_2(\mathbb{C})$ et $\mbox{SO}_3(\mathbb{R})$ : un sous-groupe et ses conjugués, d'indice fini dans leur normalisateur et d'intersection mutuelle triviale, recouvrent ``génériquement'' le groupe ambiant. Dans la perspective de la théorie des modèles, nous travaillons avec l'hypothèse de l'existence d'une bonne notion dimension sur les ensembles définissables (il faut distinguer le cas o-minimal et le cas rangé). Nous accordons une attention particulière au cas rangé. En étudiant la géométrie d'incidence induite par les involutions, nous esquissons une classification des quasi groupes de Frobenius et nous déterminons ainsi sous quelles conditions des groupes classiques peuvent être identifiés dans un cadre dimensionnel.

    \begin{center}
        \textbf{Abstract}
    \end{center}
   We are interested in a class of groups, quasi Frobenius groups (with involutions), whose internal structure generalizes that of the classical groups $\mbox{GA}_1(\mathbb{C})$, $\mbox{PGL }_2(\mathbb{C})$ and $\mbox{SO}_3(\mathbb{R})$: a subgroup and its conjugates, of finite index in their normalizer and trivial mutual intersection, cover ``generically'' the ambient group. From the perspective of model theory, we work with the hypothesis of the existence of a good notion of dimension on definable sets (we must distinguish between the o-minimal case and the ranked case). We pay special attention to the ranked case. By studying the geometry of incidence induced by involutions, we sketch a classification of quasi Frobenius groups and thus determine under which conditions classical groups can be identified in a dimensional framework.
\end{abstract}

\section{Introduction}
Les groupes algébriques sur un corps algébriquement clos sont des exemples importants de groupes rangés (aussi appelés groupes de rang de Morley fini); leurs ensembles définissables peuvent être munis d'une dimension finie aux propriétés remarquables. \` A cet égard, il est possible d'établir une analogie avec les groupes de Lie réels et plus généralement avec les groupes définissables dans les structures o-minimales. Ces derniers peuvent être également décrits du point de vue d'une fonction de dimension qui présente de nombreux points communs avec le rang de Morley (elle est additive et définissable). Dans les deux cas, nous disposons d'une notion appropriée de généricité qui s'applique aux parties définissables de dimension maximale. Dans la perspective de la théorie des modèles, il semble donc possible de décrire dans un cadre commun - celui des groupes \emph{dimensionnels}, dont la définition est donnée en §2 - les groupes de Lie réels et les groupes algébriques sur un corps algébriquement clos. Comme références générales pour les groupes rangés, nous renvoyons à \cite{BN} et \cite{ABC08}; en ce qui concerne les groupes définissables dans une structure o-minimale, on pourra consulter \cite{Ote}.

On peut illuster ce principe méthodologique en considérant les groupes $\mbox{SO}_3(\mathbb{R})$ et $\mbox{PGL}_2(\mathbb{C})$. Ce sont des groupes dimensionnels qui présentent une structure interne très proche, se reflétant dans la géométrie d'incidence induite par les involutions. Comme le font remarquer J. Wiscons et A. Deloro dans \cite{DW} et avant eux A. Nesin, les involutions forment un plan projectif dans $\mbox{SO}_3(\mathbb{R})$; les axiomes sont satisfaits seulement génériquement dans $\mbox{PGL}_2(\mathbb{C})$. Il s'agira de déterminer dans quelle mesure cette géométrie des involutions permet de caractériser ces groupes classiques parmi les groupes dimensionnels. 

Pour établir cette théorie commune, nous sommes conduits à introduire les définitions suivantes :
\begin{definition}\label{def quasi Frobenius}
Un groupe $G$ contenant un sous-groupe propre $C$ est un quasi groupe de Frobenius si $C$ est d'indice fini dans son normalisateur et TI, \textit{i.e.}, $C^g\cap C=\{ 1 \}$ pour tout $g\notin N_G(C)$. On dit que $C$ est un quasi complément de Frobenius. 
\end{definition}
\begin{definition}
Soit $C<G$ un quasi groupe de Frobenius. On note $n$ l'indice de $C$ dans $N_G(C)$. On dira que $G$ est :

\begin{itemize}
    \item de degré pair si $n$ est un entier pair,
    \item de degré impair si $n$ est un entier impair.
   \end{itemize}
Si le degré est égal à un, on parlera plutôt de groupe de Frobenius.
\end{definition}
Dans la suite de l'article, sauf mention explicite du contraire, le terme de ``degré ''  ne désignera pas le degré de Morley. 

Nous cherchons à classifier les quasi groupes de Frobenius dimensionnels; nous insisterons en particulier sur les configurations de degré pair dont l'étude forme le c\oe ur de l'article. En général, nous supposerons que les groupes en question sont $U_2^{\perp}$ : il n'y a pas de 2-groupe abélien élémentaire infini. On dira qu'un groupe contenant des involutions et $U_2^{\perp}$ est de \textit{type impair}. Du point de vue des quasi groupes de Frobenius de degré pair, nous pouvons formuler les deux conjectures suivantes :
\begin{conjecture}
\begin{enumerate}
\item \cite{DW} 
Soit $C<G$ un quasi groupe de Frobenius rangé connexe tel que $C$ est définissable et connexe. Supposons que $G$ est de degré pair et de type impair. Alors $G$ est isomorphe à $\mbox{PGL}_2(K)$ pour un corps algébriquement clos $K$ de caractéristique différente de deux.
\item Soit $C<G$ un quasi groupe de Frobenius définissable dans une structure o-minimale, connexe, tel que $C$ est définissable et connexe. Supposons que $G$ contient des involutions. Si $\bigcup_G C^g$ contient toutes les translations, \textit{i.e.}, les produits de deux involutions, alors $G$ est isomorphe à $\mbox{SO}_3(R)$ pour un corps réel clos $R$.
\end{enumerate}
\end{conjecture}
L'étude de ce type de configurations trouve son origine dans les articles \cite{Nes} et \cite{NPR}. Les groupes $\mbox{PGL}_2(K)$ et $\mbox{SO}_3(R)$ y sont identifiés parmi les groupes dimensionnels de petite dimension (inférieure ou égale à trois). 

Nous évoquerons également les groupes de Frobenius mais nous restreignons notre attention au cas rangé. La conjecture suivante guidera notre étude :
\begin{conjecture}\cite{BN}
Un groupe de Frobenius rangé $C<G$ est scindé, \textit{i.e.}, s'écrit $U\rtimes C$ pour un sous-groupe définissable $U$.
\end{conjecture}
Cette conjecture est vérifiée pour les groupes de Frobenius algébriques et finis.
Dans le cas rangé, les travaux de A. Nesin constituent une première élucidation de la structure des groupes de Frobenius. Remarquons que $\mbox{GA}_1(\mathbb{C})\simeq \mathbb{C}^{+}\rtimes \mathbb{C}^{\times}$, le groupe des transformations affines complexes, est un exemple de groupe de Frobenius rangé connexe. Si les involutions induisent également une géométrie d'incidence, cette dernière se distingue assez nettement de la géométrie observée précédemment. Pour plus de détails, nous renvoyons à \cite{CT}. En ce qui concerne les structures o-minimales, la classification des groupes de Frobenius et en particulier des groupes strictement 2-transitifs est clarifiée par les articles \cite{Tent} et \cite{MMT}. En effet, dans \cite{Tent}, il est établi que les groupes strictement 2-transitifs définissables dans une struture o-minimale sont isomorphes aux groupes des transformations affines d'un corps gauche (interprétable); ce sont en particulier des groupes de Frobenius scindés.

Dans cet article, nous introduisons d'abord le cadre approprié pour définir la notion de groupe dimensionnel et nous montrons sous quelles conditions il est possible d'exploiter la structure des quasi groupes de Frobenius pour interpréter $SO_3(K)$, pour un corps $K$ interprétable (§2).

Nous étudions ensuite plus précisément le cas o-minimal et nous démontrons le théorème suivant (§3) :
\begin{theoreme}\label{classification QFP o-minimal}
Soit $C<G$ un quasi groupe de Frobenius définissable dans une structure o-minimale, connexe, tel que $C$ est définissable et connexe. Supposons que $G$ contient des involutions. On suppose que $\bigcup_G C^g$ contient tous les 2-éléments et toutes les translations (produits de deux involutions). Alors : 
\begin{enumerate}
    \item $G$ est un groupe semi-simple, \textit{i.e.}, ne contient pas de sous-groupe abélien normal infini.
    \item Si on suppose de plus $C$ nilpotent, alors $G/Z(G)\simeq \mbox{SO}_3(R)$ pour un corps réel clos $R$ interprétable. 
\end{enumerate}
\end{theoreme}
Nous développons ensuite notre analyse dans le contexte rangé (§4). Nous formulons d'abord quelques résultats généraux sur les quasi groupes de Frobenius (structure  induite, structure de la 2-torsion, résolubilité et étude du groupe de Weyl). Une attention particulière est ensuite portée à l'étude des translations. Nous montrons sous l'hypothèse de la générosité des sous-groupes de Borel (l'ensemble  des éléments formé par un sous-groupe définissable connexe résoluble  maximal et ses conjugués est générique dans le groupe ambiant) les deux théorèmes suivants :
\begin{theoreme}\label{classification QFP}
Soit $C<G$ un quasi groupe de Frobenius rangé connexe tel que $C$ est définissable. Supposons que $G$ est de degré pair et de type impair. Si les sous-groupes de Borel sont généreux, alors $G\simeq \mbox{PGL}_2(K)$ pour un corps $K$ (rangé) algébriquement clos de caractéristique différente de 2.
\end{theoreme}
\begin{theoreme}\label{classification frobenius}
Soit $C<G$ un groupe de Frobenius rangé connexe de type impair avec $C$ résoluble et $G$ non résoluble. Si les sous-groupes de Borel sont généreux, alors $C$ n'est pas nilpotent.
\end{theoreme}
\` A cet égard, soulignons le parallélisme entre notre approche et l'étude des configurations de type ``CiBo'' (voir \cite{CJ}, \cite{Del} et \cite{DJ3}). Les théorèmes précédents constituent un analogue de certains résultats de \cite{CJ} et \cite{DJ3} où une hypothèse de minimalité (``minimal simple connexe'' ou ``$N^{\circ}_{\circ}$'') a été remplacée par l'hypothèse ``quasi groupe de Frobenius''. Les techniques employées dans ces travaux jouent un rôle prépondérant dans notre classification des quasi groupes de Frobenius rangés.

Cet article se base sur un chapitre de la thèse doctorale de l'auteur, effectuée sous la direction de Frank Wagner. Nous souhaiterions le remercier chaleureusement pour ses relectures attentives et ses nombreuses remarques qui se sont avérées essentielles en bien des points. Nous remercions également le relecteur anonyme pour ses nombreux commentaires qui ont permis de clarifier certaines démonstrations et d'améliorer la présentation et l'organisation de l'article.
\section{Groupes dimensionnels et théorème de Bachmann}
\'Etant donné une structure $M$, on dira qu'un groupe $G$ définissable dans $M$ est \textit{dimensionnel} si on peut munir les ensembles définissables de $G$ et de ses puissances cartésiennes d'une dimension finie $dim$ satisfaisant les axiomes suivants (voir \cite{EJO1} et \cite{EJO2}) :
\begin{enumerate} 
    \item[($A_1$)] Si $f$ est une fonction définissable entre deux sous-ensembles définissables $A$ et $B$, alors $\{ b\in B : dim(f^{-1}(b))=n \}$ est définissable.
    \item[($A_2$)] Si $f$ est une fonction définissable surjective entre deux sous-ensembles définissables $A$ et $B$, telle que les fibres sont de dimension constante égale à $m$, alors $dim(A)=dim(B)+m$.
    \item[($A_3$)] $A$ est fini si et seulement si $dim(A)=0$.
    \item[($A_4$)] $dim(A\cup B)=max\{dim(A),dim(B)\}$.
\end{enumerate}
Les groupes de rang de Morley fini, les groupes définissables dans une structure o-minimale, dans le corps des nombres p-adiques, dans un corps P-minimal sont des exemples de groupes dimensionnels.

On parlera de groupe $dimensionnel^{+}$ si de plus le groupe satisfait la condition de chaîne descendante pour les sous-groupes définissables et si la dimension peut être étendue aux sous-ensembles interprétables du groupe. En particulier, on peut attacher une dimension aux sections définissables (les quotients d'un sous-groupe définissable par un autre sous-groupe définissable).

Les groupes définissables dans une structure o-minimale \cite[Théorème 7.2]{ED} et les groupes rangés (par définition) sont des exemples de groupes $dimensionnels^{+}$ (mais ce n'est plus le cas des groupes définissables dans le corps des nombres p-adiques, par exemple).

Un groupe $dimensionnel^{+}$ $G$ admet une \textit{composante connexe} $G^{\circ}$, \textit{i.e.}, un plus petit sous-groupe définissable d'indice fini : il suffit de considérer l'intersection de tous les sous-groupes définissables d'indice fini. On dira qu'il est \textit{connexe} si $G^{\circ}=G$. De plus, chaque sous-ensemble $X\subseteq G$ est contenu dans un plus petit sous-groupe définissable, appelé \textit{enveloppe définissable} et noté $d(X)$. Remarquons que, comme dans le cas rangé, certaines propriétés sont préservées lorsque l'on passe à la clôture définissable d'un sous-groupe (non-définissable) $H$ \cite[Lemmes 3.3 et 3.4]{EJO2} :
\begin{itemize}
    \item Soit $K\subseteq G$. Si $H$ est $K$-invariant, alors $d(H)$ est $K$-invariant.
    \item Si $H$ est nilpotent de classe $n$, alors $d(H)$ est nilpotent de classe $n$.
\end{itemize}
La condition de chaîne descendante pour les sous-groupes définissables rend également possible le relèvement de la torsion \cite[Fait 3.9]{EJO1} :

Soit $H$ un sous-groupe définissable normal tel que $xH$ est un $p$-élément dans $G/H$. Alors il existe un $p$-élément dans $xH$. 

L'existence d'une dimension finie permet de caractériser facilement les parties qui sont, en un certain sens, ``larges''. On dira qu'une partie $Y$ est \textit{générique} si elle contient une partie définissable $X$ telle que $dim(X)=dim(G)$. On dira que $X$ est \textit{large} si $dim(G\setminus X)<dim(G)$. De plus, étant donné un sous-groupe définissable $H$, l'équation $dim(G)=dim(G/H)+dim(H)$ montre qu'un sous-groupe définissable $H$ est d'indice fini dans $G$ ssi $H$ est générique dans $G$.

Passons maintenant à l'étude des quasi groupes de Frobenius $dimensionnels^{+}$. La notion de quasi groupe de Frobenius a été définie à la Définition \ref{def quasi Frobenius}. Soit $C<G$ un quasi groupe de Frobenius $dimensionnel^{+}$ connexe tel que $C$ est définissable et connexe. On remarque que $N_G(C)^{\circ}=C$ et que $C_G(g)^{\circ}\leq C$, pour tout $1\neq g\in C$ : $C_G(g)\leq N_G(C)$ par la propriété TI et $C_G(g)^{\circ}\leq N_G(C)^{\circ}=C$.
\begin{nota}
Soit $C<G$ un quasi groupe de Frobenius. Pour $x\in \bigcup_G C^g$, on note $C_x$ le conjugué de $C$ contenant $x$.
\end{nota} 
\begin{Lemme}\label{genericté frobenius}
Soit $C<G$ un quasi groupe de Frobenius $dimensionnel^{+}$ connexe tel que $C$ est définissable et connexe. Alors $\bigcup_G C^g$ est générique dans $G$.
\end{Lemme}
\begin{proof}  En raisonnant comme dans \cite{Jal1}, l'additivité et la définissabilité de la dimension permettent de conclure.\end{proof}
Avant de poursuivre l'analyse des quasi groupes de Frobenius $dimensionnels^{+}$ connexes, il nous faut souligner une différence essentielle entre le contexte rangé et le contexte o-minimal. Pour un groupe connexe rangé $G$, la notion de largeur et celle de généricité coïncident alors que ce n'est pas forcément le cas pour un groupe connexe définissable dans une structure o-minimale. En particulier, deux parties génériques ne s'intersectent pas forcément non-trivialement. 

De plus, l'action des automorphismes définissables involutifs peut être caractérisée de façon commode dans le contexe rangé (mais pas forcément dans le cas o-minimal) : si un automorphisme involutif définissable d'un groupe rangé connexe a un nombre fini de points fixes, alors l'action se fait par inversion et le groupe est abélien. Or, la compréhension de l'action des involutions se révèle indispensable à notre étude des quasi groupes de Frobenius $dimensionnels^{+}$ connexes. Nous sommes donc conduits à proposer un nouvel axiome pour identifier une sous-classe des groupes $dimensionnels^{+}$:
\vspace{0.3cm}
\begin{itemize}
\item[($A_5$)] Pour tout sous-groupe définissable connexe $H$ et tout automorphisme involutif définissable $i$ de $H$, si $C_H(i)$ est fini, alors $H$ est abélien inversé par $i$.
\end{itemize}
\vspace{0.3cm}
Nous ne savons pas si un tel axiome (vérifié par les groupes rangés) reste valable pour une sous-classe naturelle de groupes définissables dans une structure o-minimale (par exemple, pour les groupes définissablement simples et/ou définissablement compacts). Dans la suite, nous indiquerons à quel moment $(A_5)$ semble requis.
\begin{Remarque}\label{remarque}
Soient $C$ un sous-groupe abélien définissable connexe et $i$ un automorphisme involutif définissable; alors $(A_5)$ n'est plus nécessaire pour montrer que l'action est par inversion. En effet, posons $X=\{ c^ic^{-1}: c\in C \}$ (c'est un ensemble définissable inversé par $i$). Les fibres de l'application définissable $\phi : C\longrightarrow X$ définie par $\phi (c)=c^ic^{-1}$, ont même dimension que $C_C(i)$, qui est un groupe fini. Par conséquent, l'additivité de la dimension nous donne $dim(X)=dim(C)$. Le sous-groupe formé par les éléments inversés par $i$ est générique et donc il est égal à $C$ par connexité.
\end{Remarque}

Nous allons poursuivre l'étude des translations menée dans \cite{DW} et nous utiliserons les méthodes de l'école de Bachmann pour interpréter $\mbox{SO}_3(K)$ dans un cadre dimensionnel.

\textit{Jusqu'à la fin de la section, soit $C<G$ un quasi groupe de Frobenius $dimensionnel^{+}$ connexe tel que $C$ est définissable et connexe. Supposons que $G$ est de type impair et que $\bigcup_G C^g$ contient tous les 2-éléments. Dans le cas rangé, l'hypothèse $U_2^{\perp}$ suffit (voir Lemme \ref{générosité complément}). On suppose également que toutes les translations (les produits de deux involutions) sont contenues dans $\bigcup_G C^g$.}
\begin{Lemme}
Le sous-groupe $C$ est d'indice pair dans son normalisateur. De plus, il existe une involution $1\neq k\in N_G(C)\setminus C$.
\end{Lemme}
\begin{proof}
Soit $i$ une involution contenue dans $C$ et $j\neq i$ une autre involution contenue dans un conjugué $C_j$ distinct de $C$. Par hypothèse, la translation $x=ij$ est contenue dans un conjugué $C_x$ de $C$. On a donc $C_G(x)^{\circ}\leq C_x$; de plus, $i$ et $j$ appartiennent à $N_G(C_G(x)^{\circ})\leq N_G(C_x)$. Si $i$ et $j$ étaient contenues dans $C_x$, alors on aurait d'une part $C_x\cap C\neq \{ 1 \}$ et d'autre part $C_x\cap C_j\neq \{ 1 \}$ et finalement $C=C_x=C_j$, contradiction.
Par conséquent, $i$ et $j$ appartiennent à $N_G(C_x)\setminus C_x$.
\end{proof}
On dit qu'un groupe est \textit{semi-simple} s'il n'a pas de sous-groupe abélien normal infini.
Soit $C<G$ un quasi groupe de Frobenius $dimensionnel^{+}$ connexe tel que $C$ est définissable et connexe. Alors $G$ a un centre fini. En effet, $Z(G)\leq N_G(C)$, donc $Z(G)^{\circ}\leq C$. En prenant des conjugués, $Z(G)^{\circ}=\{ 1 \}$. 
\begin{Lemme}\label{semi-simple}
 Le groupe $G$ est semi-simple.
\end{Lemme}
\begin{proof}
Soit $A$ un sous-groupe abélien normal infini. Quitte à considérer la composante connexe de son enveloppe définissable, on peut supposer que $A$ est définissable et connexe. Si $A\cap C\neq \{ 1 \}$, soit $1\neq x\in A\cap C$; alors $A\leq C_G(x)^{\circ}\leq C$ et donc $C=G$, contradiction. On a donc $A\cap \bigcup_G C^g=\{ 1 \}$. Le groupe $A$ est en particulier 2-divisible car il ne contient pas d'involutions. 
Soit $i$ une involution; comme on peut supposer que $i\in C$, $C_A(i)^{\circ}\leq A\cap C=\{1\}$, le sous-groupe $C_A(i)$ est donc fini et il existe un élément $a\in A$ tel que $a^ia^{-1}=b\in A$ est non-trivial et inversé par $i$. Mais alors $1\neq b^2\in A$ et $b^2=bib^{-1}i=i^{b^{-1}}i\in \bigcup_{g\in G} C^g$, contradiction.
Par conséquent, il n'y a pas de sous-groupe abélien normal infini.
\end{proof}
On suppose désormais que $G$ satisfait également $(A_5)$.
\begin{Lemme}\label{inversion}
Soit $k$ une involution telle que $k\in N_G(C)\setminus C$. Alors $C_C(k)$ est fini. De plus, $C$ est abélien 2-divisible inversé par $k$ (et donc $C\subseteq I\cdot I$).
\end{Lemme}
\begin{proof}
 Supposons que $C_C(k)$ est infini, on a donc $C_C(k)^{\circ}\leq C_k\cap C$ et finalement $C_k=C$, contradiction. Par $(A_5)$, l'action de $k$ est donc bien par inversion et $C$ est abélien. C'est un groupe 2-divisible car il contient un nombre fini d'involutions : l'élévation au carré définit un homomorphisme de groupe définissable de noyau fini. Puisque le groupe est connexe, c'est un homomorphisme surjectif.
\end{proof}
\begin{Remarque}
Nous n'utilisons plus l'axiome $(A_5)$ dans la suite.
\end{Remarque}
\begin{Lemme}
Le sous-groupe (abélien) $C$ contient une unique involution et toutes les involutions de $N_G(C)\setminus C$ appartiennent au même translaté de $C$.
\end{Lemme}
\begin{proof}
Supposons par l'absurde que $i\neq j$ sont deux involutions distinctes appartenant à $C$.  Soit $\alpha \in N_G(C)\setminus C$ une involution. Puisque l'involution $\alpha$ agit par inversion, elle commute avec l'involution $x=ij$; mais $i$ appartient à $N_G(C_{\alpha})\setminus C_{\alpha}$ donc inverse $C_\alpha$, $j$ de même, donc $x\in C_G(C_{\alpha})$. D'après le Lemme \ref{inversion}, $x\in C_{\alpha}\cap C=\{1\}$, contradiction. Ceci montre l'unicité.
 
Soient $j,k\in N_G(C)\setminus C$; puisque l'action est par inversion, $jk\in C_G(C)$. Mais alors $C\leq C_G(jk)^{\circ}\leq C_{jk}$ et donc $C_{jk}=C$ et $jk\in C$. Finalement, $k\in jC$ et donc toutes les involutions de $N_G(C)\setminus C$ sont égales modulo $C$. Réciproquement, le translaté $jC$ est composé uniquement d'involutions.
\end{proof}
En procédant comme dans la démonstration du Théorème 3.31 de \cite{Poi1}, on peut munir l'ensemble des involutions d'une structure de plan projectif : à chaque involution $i$ (point), on associe la ligne $D_i=\{ k\in I : ki=ik\}\setminus \{i\}=\{k\in I : k\in N_G(C_i)\setminus C_i \}=kC_i$ pour chaque $k\in N_G(C_i)\setminus C_i$. 

Dans ce cas, par deux points distincts $i\neq j$ passe une unique ligne à savoir $D_k$ où $k$ est l'unique involution de $C_{ij}$. En effet, $i$ et $j$ inversent $C_{ij}=C_k$ mais n'appartiennent pas à $C_{ij}$. De plus, si $i,j\in D_{\ell}\neq D_{k}$, alors $i,j\in N_G(C_{\ell})\setminus C_{\ell}$ et donc $x=ij\in C_{\ell}$, finalement $C_{\ell}=C_x=C_k$, contradiction.

De plus, deux lignes $D_i$ et $D_j$ s'intersectent en un seul point $k$, \textit{i.e.}, l'unique involution contenue dans $C_{ij}$. En effet, pour cette involution, on a bien $k\in N_G(C_i)\setminus C_i$ et $k\in N_G(C_j)\setminus C_j$ car $k\in C_G(i)\cap C_G(j)$. Soit $\ell \neq k\in D_j\cap D_i$; alors $\ell k\in C_i\cap C_j$, et $C_i=C_j$, contradiction. En particulier, pour deux involutions distinctes $i\neq j$, il existe une unique involution $k$ distincte de $i$ et de $j$ qui commute avec les deux. 

Finalement, on remarque que trois involutions distinctes $i,j,k$ sont colinéaires ssi $ijk$ est une involution : $ijk$ est une involution [ssi $ij$ est inversé par $k$ et $ij\neq k$] ssi $k,i,j\in D(\ell)$ où $\ell$ est l'involution contenue dans $C_{ij}$. On peut dès lors appliquer une forme du théorème de Bachmann : 
\begin{fait}\cite[Fait 8.15]{BN},\cite{Schr}
Soit $G$ un groupe contenant des involutions. Si on peut munir l'ensemble des involutions $I$ d'une structure de plan projectif tel que trois involutions $i,j$ et $k$ sont colinéaires ssi $ijk$ est une involution, alors le sous-groupe $\langle I \rangle=I\cdot I$ est isomorphe à $\mbox{SO}_3(K,f)$ pour un corps $K$ interprétable dans $G$ et une forme quadratique anisotrope $f$ définie sur $K^3$.
\end{fait}
Il suffit de remarquer que $\langle I \rangle = I\cdot I=\bigcup_G C^g$ est un sous-groupe définissable de même dimension que $G$ (Lemme \ref{genericté frobenius}) et donc lui est égal par connexité. Notamment, $G\simeq \mbox{SO}_3(K,f)$.
 
 Dans un univers rangé, un corps interprétable K est algébriquement clos \cite{Mac} et toute forme quadratique définie sur un K-espace vectoriel de dimension supérieure ou égale à deux est alors forcément isotrope. Par conséquent, comme corollaire, on obtient le Théorème A de \cite{DW} :
 \begin{theoreme}\label{transpositions QFP}
 Soit $C<G$ un quasi groupe de Frobenius rangé connexe tel que $C$ est définissable et connexe. Supposons que $G$ est de type impair. Alors $\bigcup_G C^g$ ne contient pas toutes les translations.
 \end{theoreme}

\section{Le cas o-minimal} 
Cette section est consacrée à la démonstration du Théorème \ref{classification QFP o-minimal} :

Soit $C<G$ un quasi groupe de Frobenius définissable dans une structure o-minimale, connexe, tel que $C$ est définissable et connexe. Supposons que $G$ contient des involutions. On suppose que $\bigcup_G C^g$ contient tous les 2-éléments et toutes les translations (produits de deux involutions). Alors : \begin{enumerate}
\item $G$ est un groupe semi-simple, \textit{i.e.}, ne contient pas de sous-groupe abélien normal infini.
\item Si on suppose de plus $C$ nilpotent, alors $G/Z(G)\simeq \mbox{SO}_3(R)$ pour un corps réel clos $R$ interprétable. 
\end{enumerate}
Notons que nous affaiblissons l'hypothèse de commutativité de la Remarque  \ref{remarque} en supposant que le quasi complément est seulement nilpotent.

Nous renvoyons à \cite{PPS1} pour les théorèmes concernant la structure des groupes semi-simples définissables dans une structure o-minimale. Nous mobilisons également la théorie des sous-groupes de Carter et de Cartan développée dans \cite{EJO2}. Un \textit{sous-groupe de Cartan} $Q$ est un sous-groupe nilpotent maximal tel que pour tout sous-groupe $X\leq Q$ normal d'indice fini, $X$ est un sous-groupe d'indice fini dans $N_G(X)$. D'après \cite{EJO2}, on sait que les sous-groupes de Cartan d'un groupe définissable dans une structure o-minimal sont définissables. De plus, les \textit{sous-groupes de Carter}, \textit{i.e.}, les sous-groupes définissables connexes nilpotents d'indice fini dans leur normalisateur correspondent aux composantes connexes des sous-groupes de Cartan. Il est possible de caractériser précisément la structure des sous-groupes de Cartan des groupes définissablement simples définissables dans une structure o-minimale (et plus généralement des sous-groupes définissables semi-simples).
\begin{fait}\label{Cartan} \cite[Théorème 57]{EJO2}
Soit $G$ un groupe définissablement simple définissable dans une structure o-minimale. Alors les sous-groupes de Cartan sont abéliens et ont tous même dimension.
\end{fait}
\textit{Soit $C<G$ un quasi groupe de Frobenius définissable dans une structure o-minimale, connexe, tel que $C$ est définissable et connexe. Supposons que $G$ contient des involutions et que $\bigcup_G C^g$ contient les 2-éléments et les translations. Nous savons déjà que $G$ est semi-simple (Lemme \ref{semi-simple}). Nous supposons de plus que $C$ est nilpotent (c'est donc un sous-groupe de Carter).}
\begin{Lemme}
Le groupe $G/Z(G)$ est un quasi groupe de Frobenius avec $CZ(G)/Z(G)$ comme quasi complément de Frobenius.
\end{Lemme}
\begin{proof}
On note que $Z(G)$ est fini. Supposons $\alpha\in \overline{C}\cap \overline{C^g}\setminus\{\overline{1}\}$, relevé par $a$. Alors $a=c_1z_1=c_2^gz_2$ avec $c_1, c_2\in C\setminus\{1 \}$, et $z_1, z_2\in Z(G)$. On a $C_G(a)^{\circ}=C_G(c_1)^{\circ}\leq C$ et donc $C=C^g$; ceci montre que $\overline{C}$ est TI. Soit maintenant $\nu \in N_{\overline{G}}(\overline{C})$, relevé par $n\in G$. Alors $C^n\leq C\cdot Z(G)$, donc par connexité $C^n\leq C$ et $n\in N_G(C)$. Le groupe $N_{\overline{G}}(\overline{C})/\overline{C}$ est donc fini.
\end{proof}
Nous pouvons passer à la démonstration du Théorème \ref{classification QFP o-minimal}. D'après \cite[Théorème 4.1]{PPS1}, le groupe $G/Z(G)$ est isomorphe à un produit fini de groupes simples définissables, les $H_i$. Comme $G/Z(G)$ a des involutions, on peut supposer que $H_1$ contient des involutions. Soit $i\in H_1\cap \overline{C}_i$, alors $H_2\leq C_{\overline{G}}(i)^{\circ}\leq \overline{C}_i$, contradiction. Le groupe $G/Z(G)$ est donc définissablement simple. D'après le Fait \ref{Cartan}, ses sous-groupes de Carter sont abéliens (et de même dimension). Or, $\overline{C}$ est un sous-groupe de Carter, c'est donc un groupe abélien. On peut désormais appliquer les résultats de la section précédente et par conséquent, $G/Z(G)\simeq \mbox{SO}_3(K,f)$.

Montrons désormais que le corps $K$ interprétable est réel clos et que la forme quadratique anisotrope correspond à un produit scalaire. Nous avons besoin du fait suivant :
\begin{fait}\cite[Théorème 3.9]{Pill88} Un corps commutatif infini $K$ définissable dans une structure o-minimale est soit réel clos (et $dim(K)=1$) soit algébriquement clos. 
\end{fait}
Puisque la forme quadratique est anisotrope, le corps $K$ est réel clos. Finalement, le théorème d'inertie de Sylvester nous permet de conclure que 
$G/Z(G)\simeq \mbox{SO}_3(K)$ où $K$ est réel clos.
\section{Le cas rangé}
\subsection{Préliminaires généraux}
 Nous rappelons désormais plusieurs faits plus ou moins spécialisés sur les groupes rangés qui seront utilisés dans cet article. On pourra consulter \cite{BN} pour une caractérisation axiomatique de la dimension dans ce contexte. Lorsqu'on parlera de groupe algébrique, on entend groupe des K-points dans un corps algébriquement clos $K$. Comme pour les structures o-minimales, nous dirons qu'une partie définissable $X$ d'un groupe rangé $G$ est générique si $RM(X)=RM(G)$.
\subsubsection{Générosité}
Rappelons que dans un groupe connexe rangé, deux parties génériques s'intersectent toujours non-trivialement. Introduisons un peu de terminologie. On dira qu'un sous-groupe définissable $H$ est \textit{généreux} si $\bigcup_{g\in G} H^{g}$ est générique. On dira également qu'un sous-groupe définissable $H$ est \textit{presque-autonormalisant} si $H$ est d'indice fini dans son normalisateur. Enfin, un sous-groupe définissable $H$ est \textit{génériquement disjoint de ses conjugués} dans $G$ si $H\setminus \bigcup_{g\in G\setminus N(H)}(H\cap H^g)$ est générique dans $H$. On rappelle plusieurs résultats concernant la générosité et la généricité (voir \cite{Jal1}).
\begin{fait}\label{générosité}
\begin{enumerate}
\item \cite[Lemme IV.1.2]{ABC08} Soient $G$ un groupe rangé  connexe et $H$ un sous-groupe définissable presque-autonormalisant et génériquement disjoint de ses conjugués. Alors $H$ est généreux dans $G$.
\item \cite[Lemme 2.4]{Jal1} Si $A\leq B \leq C$ sont trois sous-groupes définissables avec $B$ connexe, alors la générosité est transitive.
\end{enumerate}
\end{fait}
\subsubsection{Indécomposables}
Nous supposerons connus les raisonnements par indécomposable \cite[§5.4]{BN}; nous ferons également usage du résultat de Zilber concernant l'interprétation d'un corps à partir d'une action définissable sur un groupe abélien \cite[Théorème 9.1]{BN}.
\subsubsection{Théorie de Sylow}
Les groupes rangés peuvent être analysés \textit{via} leurs éléments de torsion. Il existe une théorie des $p$-sous-groupes de Sylow ($p$-sous-groupes maximaux) mais elle n'est pleinement développée que pour $p=2$ ou dans le contexte des groupes (connexes) résolubles rangés.
D'après \cite[§6.4]{BN}, les 2-Sylow sont conjugués et leurs composantes connexes s'écrivent comme produit central d'un $2$-groupe abélien divisible et d'un groupe d'exposant fini définissable connexe.
\begin{definition}
Soit $G$ un groupe rangé et soit $S^{\circ}=T*B$ la composante connexe d'un 2-sous-groupe de Sylow, où $T$ est abélien divisible et $B$ d'exposant fini.
\begin{enumerate}
    \item Si $S^{\circ}=\{1\}$, alors $G$ est de type dégénéré.
    \item Si $T\neq \{1\}$ et $B\neq \{1\}$, alors $G$ est de type mixte.
    \item Si $T\neq \{1\}$ et $B=\{1\}$, alors $G$ est de type impair.
    \item Si $T=\{1\}$ et $B\neq \{1\}$, alors $G$ est de type pair.
\end{enumerate} 
\end{definition}
Dans cette perspective, un \textit{$p$-groupe unipotent} est un $p$-groupe nilpotent définissable connexe d'exposant fini. On dira qu'un groupe rangé est $U_{p}^{\perp}$ s'il ne contient pas de $p$-groupe unipotent. Un \textit{$p$-tore} est \textit{a contrario} un $p$-groupe abélien divisible isomorphe à une somme finie de $p$-groupes de Prüfer. Le $p$-rang de Prüfer d'un groupe rangé désigne alors le $p$-rang d'un $p$-tore maximal. Le cas $p=2$ revêt une importance particulière en raison du résultat de conjugaison évoqué précédemment. De plus, bien qu'on ne dispose pas d'un analogue du théorème de Feit-Thompson qui éliminerait le cas dégénéré, les involutions, comme dans le cas des groupes simples finis, sont souvent au c\oe ur de l'analyse. Nous renvoyons à \cite[§10]{BN} pour plus de détails sur les propriétés des involutions dans notre contexte. Nous supposerons également connues les propriétés des sous-groupes définissables \textit{fortement inclus} (voir \cite[Théorème 10.19]{BN}).

Dans le cas d'un groupe connexe, la structure de la $p$-torsion satisfait la propriété additionnelle suivante :
\begin{fait}\label{$p$-Sylow infini} \cite{BBC}
Soit $G$ un groupe rangé connexe dont les $p$-sous-groupes de Sylow sont finis. Alors $G$ ne contient pas de $p$-éléments.
\end{fait}
On cite une autre conséquence de l'analyse de la torsion dans les groupes connexes :
\begin{fait}\label{centralisateur infini} \cite{BBC}, \cite[IV.Corollaire 4.18] {ABC08}
Soit $G$ un groupe rangé connexe. Alors  $C_G(g)^{\circ}$ est infini pour tout $g\in G$.
\end{fait}
\subsubsection{Groupes résolubles}
L'étude d'un groupe rangé passe en général par une compréhension de ses \textit{sous-groupes de Borel}, \textit{i.e.}, de ses sous-groupes définissables connexes résolubles maximaux. Plus généralement, nous allons nous intéresser aux propriétés des sous-groupes définissables (connexes) résolubles.
\begin{fait}\label{G' nilpotent} \cite[Corollaire 9.9]{BN}
Soit $G$ un groupe rangé connexe résoluble. Alors $G'$ est un sous-groupe définissable connexe nilpotent.
\end{fait}
\subsubsection{Notions de tore}
En exploitant l'analogie avec les groupes algébriques affines, on peut introduire pour les groupes rangés une notion de tore (éléments semi-simples) et de radical unipotent. Idéalement un groupe connexe rangé résoluble devrait pouvoir s'écrire comme un produit semi-direct d'un radical unipotent et d'un tore.

On caractérise maintenant la semi-simplicité dans le contexte des groupes rangés. Nous pourrons également mieux comprendre la nature des éléments génériques.
Un sous-groupe abélien divisible est un \textit{tore décent} si c'est la clôture définissable de sa torsion; c'est un \textit{bon tore} si cette propriété est héréditaire pour les sous-groupes définissables infinis. Le groupe multiplicatif d'un corps de caractéristique positive est un exemple typique de bon tore comme le prouve le fait suivant qui se base sur des travaux de F.O. Wagner:
\begin{fait}\label{bon tore} \cite[Proposition 4.20]{ABC08}
Soit $K$ un corps rangé de caractéristique positive. Alors $K^{\times}$ est un bon tore.
\end{fait}
Remarquons également que le résultat suivant de F.O. Wagner permet de conjecturer la non-existence des \textit{mauvais corps} en caractéristique positive :
\begin{fait}\label{p-Mersennes} \cite{Wag4}
S'il existe une infinité de $p$-nombres premiers de Mersenne, \textit{i.e.}, de la forme $(p^n-1)/(p-1)$, alors il n'existe pas de mauvais corps de caractéristique positive $p>0$.
\end{fait}
Grâce aux travaux de G. Cherlin et T. Alt\i nel puis d'O. Frécon notamment, on sait que le comportement des tores décents présente des similitudes avec celui des tores algébriques.
\begin{fait}\label{tore décent}
\begin{enumerate}
\item \cite{Che} Soit $G$ un groupe rangé. Alors les tores décents maximaux sont conjugués.
\item \cite[Lemme 3.1]{Fre3} Soient $G$ un groupe rangé et $N$ un sous-groupe normal définissable. Soit $T$ un tore décent maximal. Alors $TN/N$ est un tore décent maximal de $G/N$ et tous les tores décents maximaux ont cette forme.
\item \cite[Théorème 1]{AB} Soit $T$ un tore décent d'un groupe rangé connexe $G$. Alors $C_G(T)$ est connexe, généreux dans $G$ et $N_G(C_G(T))^{\circ}=C_G(T)$.
\end{enumerate}
\end{fait}
On cite également un fait qui décrit la $p$-torsion lorsque le groupe $U_p^{\perp}$.
\begin{fait}\label{torsion semi-simple}\cite[Théorème 3]{BC} Soit $G$ un groupe rangé connexe, et soit $g$ un $p$-élément tel que $C_G(g)^{\circ}$ est $U_p^{\perp}$. Alors $g$ appartient à un $p$-tore. 
\end{fait}
Il s'agit d'un principe de \textit{toralité}.
\subsection{Résultats généraux}
Nous commençons à élaborer une théorie générale des quasi groupes de Frobenius rangés connexes. Un certain nombre de résultats valables pour les groupes de Frobenius rangés (étudiés dans \cite{BN}) se généralisent à ce cadre plus général. Remarquons que nous faisons usage de manière cruciale d'arguments impliquant l'étude de la généricité et de la générosité dans les groupes connexes (ce qui les rend difficilement généralisables au cas o-minimal).
\subsubsection{Structure induite}
Dans un quasi groupe de Frobenius rangé connexe, l'intersection d'un sous-groupe définissable connexe avec un quasi complément de Frobenius définissable connexe induit une structure de quasi groupe de Frobenius connexe.
\begin{Lemme}\label{conjugaison complément Frobenius}(à comparer avec \cite[Lemme 11.10]{BN})
Soit $C<G$ un quasi groupe de Frobenius rangé connexe tel que $C$ est définissable et connexe. Soit $H<G$ un sous-groupe définissable connexe. Si $1<(H\cap C)<H$ alors $(H\cap C)<H$ est encore un quasi groupe de Frobenius. De plus, $H\cap C$ est infini et connexe. Finalement, pour tout conjugué $C_1$ de $C$ tel que $(H\cap C_1)\neq \{ 1 \}$, il existe $1\neq h\in H$, tel que $H\cap C_1=(H\cap C)^h$.
\end{Lemme}
\begin{proof}
Montrons que $H\cap C$ est TI dans $H$. Supposons que $(H\cap C)^h\cap (H\cap C)\neq \{ 1 \}$ pour un $h\in H$; on a donc $C^h\cap C\neq \{ 1 \}$ et par conséquent $C^h=C$. Le même raisonnement montre que $N_H(H\cap C)\leq N_G(C)$ et donc $N_H(H\cap C)^{\circ}\leq N_G(C)^{\circ}=C$, \textit{i.e.},\ $N_H(H\cap C)^{\circ}\leq H\cap C$. Le groupe $H\cap C$ est donc TI dans $H$ et d'indice fini dans son normalisateur dans $H$.
On remarque que pour $1\neq c\in H\cap C$, $C_H(c)\leq N_H(H\cap C)$ et donc $C_H(c)^{\circ}\leq N_H(H\cap C)^{\circ}$. Or, d'après le Fait \ref{centralisateur infini}, $C_H(c)$ est infini donc $H\cap C$ est infini.

Si $H\cap C$ n'était pas connexe alors $(H\cap C)^{\circ}$ et $(H\cap C)\setminus (H\cap C)^{\circ}=X$ seraient généreux. En effet, $(H\cap C)^{\circ}$ et $X$ sont invariants sous l'action de $N_H(H\cap C)$; de plus, si $X^h$ ou $(H\cap C)^{{\circ}^h}$ intersecte $(H\cap C)$ non-trivialement alors $(H\cap C)\cap (H\cap C)^h$ est non-trivial et donc $h\in N_H(H\cap C)$. Par conséquent, les ensembles disjoints $\bigcup_H X^h$ et $\bigcup_H (H\cap C)^{{\circ}^h}$ sont génériques, ce qui contredit la connexité de $H$.
Soit $(H\cap C_1)\neq \{ 1 \}$. Par ce qui précède, c'est un sous-groupe généreux dans $H$, tout comme $(H\cap C)$. Par connexité de $H$, il existe $1\neq h\in H$ tel que $(H\cap C)\cap (H\cap C_1)^h\neq \{ 1 \}$. Par conséquent, $C_1^h=C$ et finalement $(H\cap C)=(H\cap C_1)^h$. 
\end{proof}
\begin{Remarque}
\` A notre connaissance, le quasi complément de Frobenius pourrait ne pas être définissable dans la pure structure de groupe, contrairement à ce qui se passe dans les groupes de Frobenius rangés \cite[Lemme 11.19]{BN} . Cependant, en raisonnant comme dans la démonstration précédente, on obtient le résultat suivant : si $C<G$ est un quasi groupe de Frobenius rangé connexe, avec $C$ définissable, alors $C$ est infini connexe. Dans la suite, nous préciserons seulement que le quasi complément de Frobenius est définissable puisque la connexité en découle.
\end{Remarque}
\subsubsection{Un critère de scission}
On démontre un critère de scission pour les quasi groupes de Frobenius rangés connexes avec complément définissable.
\begin{proposition} \label{scission quasi-frobenius}  (à comparer avec le Lemme 11.21 de \cite{BN})
Soit $C<G$ un quasi groupe de Frobenius rangé connexe tel que $C$ est définissable. Supposons qu'il existe des involutions et qu'elle sont toutes contenues dans $\bigcup_G C^g$. On suppose qu'il existe un groupe $A$ définissable connexe normal infini d'intersection triviale avec $C$. Alors il existe un groupe définissable connexe abélien $A_1$ tel que $G=A_1\rtimes C$.
\end{proposition}
\begin{proof}
Soit $i\in C$ une involution; alors $i$ agit par inversion sur $A$; sinon, $C_A(i)^{\circ}\leq C$, contradiction. Le groupe $A$ est donc abélien uniquement 2-divisible; si $A$ contenait une involution $j$, alors on aurait $A\leq C_j$, contradiction.
On considère maintenant $C_G(A)^{\circ}$. On a bien $C_G(A)^{\circ}\cap C_i=\{ 1 \}$; sinon, soit $1\neq x\in C_G(A)^{\circ}\cap C_i$; alors $A\leq C_G(x)^{\circ}\leq C_i$, contradiction. L'involution $i$ agit aussi par inversion sur $C_G(A)^{\circ}$, qui est donc abélien. On a de plus $C_G(C_G(A)^{\circ})^{\circ}=C_G(A)^{\circ}$, donc quitte à remplacer $A$ par $C_G(A)^{\circ}$, on peut supposer que $C_G(A)^{\circ}=A$.

Soit $g\in G$; alors $[i,g]\leq C_G(A)$. En effet, $i$ et $i^g$ agissent par inversion et $[i,g]=ii^g$. Par conséquent, $[i,G]\leq C_G(A)^{\circ}=A$. Soit $g\in G$, on considère $a=[i,g]$; puisque $A$ est 2-divisible, il existe $a'\in A$ tel que $a=a'^{2}=[i,a']$. Par conséquent, on a :
\[ [i,ga'^{-1}]=ia'g^{-1}iga'^{-1}=a'^i[i,g]a'^{-1}=a'^{-1}a'^2a'^{-1}=1. \]
On a donc $G=AC_G(i)=C_G(i)A$, mais $C_G(i)\cap A=\{ 1 \}$ (sinon, $A$ contiendrait un élément $a$ tel que $a^i=a=a^{-1}$, \textit{i.e.}, une involution, contradiction). Par conséquent, $G=A\rtimes C_G(i)$. Finalement, on déduit que $C_G(i)=C_G(i)^{\circ}=C$ et donc $G=A\rtimes C$.
\end{proof}
\subsubsection{2-torsion}
Ce sont les quasi groupes de Frobenius contenant des involutions qui nous intéressent principalement. Nous montrons un certain nombre de résultats qui clarifient la situation de la 2-torsion dans ce contexte.
\begin{Lemme} \label{générosité complément}  \cite[Proposition 1]{DW} Soit $C<G$ un quasi groupe de Frobenius rangé connexe tel que $C$ est définissable. Supposons que $G$ est de type impair. Alors les 2-éléments appartiennent à des tores décents et sont contenus dans $\bigcup_G C^g$.\end{Lemme}
Pour les groupes de Frobenius connexes, nous avons la situation suivante :
\begin{fait}\cite[Lemme 11.20]{BN} \label{involution Frobenius} Soit $C<G$ un groupe de Frobenius rangé connexe avec involutions. Alors les involutions de $G$ sont conjuguées et $C$ contient au plus une involution.
\end{fait}
\begin{proposition} \label{2-torsion Frobenius}
Soit $C<G$ un groupe de Frobenius rangé connexe avec involutions. Alors $G$ est $U_2^{\perp}$ ssi $C$ contient des involutions.
\end{proposition} 
\begin{proof}
Supposons que $G$ contient un 2-groupe abélien élémentaire infini $A$. Supposons par l'absurde que $C$ contient des involutions. Puisque par le Fait \ref{involution Frobenius} les involutions forment une classe de conjugaison, on a $I\subseteq \bigcup _G C^g$. En particulier, il existe $1\neq x\in A\cap C$ (quitte à conjuguer) et $A\leq C_G(x)\leq C$. Mais $C$ contient au plus une involution par le Fait \ref{involution Frobenius}, contradiction. Donc $C$ ne contient pas d'involutions.

Réciproquement, supposons que $C$ ne contient pas d'involutions. Alors les 2-sous-groupes de $Sylow^{\circ}$ de $G$ sont d'exposant fini (et en particulier contiennent des 2-groupes abéliens élémentaires infinis).  En effet, soit $S$ un 2-sous-groupe de $Sylow^{\circ}$ et supposons qu'il contient un 2-tore; on peut l'étendre en un tore décent maximal $T$. Mais $C_G(T)$ est généreux et il existe donc $1\neq y\in C_G(T)\cap C^g$. Par conséquent, $T\leq C_G(y)\leq C^g$, contradiction.\end{proof}
\begin{corollaire} Soit $G$ un groupe de Frobenius rangé connexe. Alors $G$ ne peut pas être de type mixte.\end{corollaire}
\begin{Remarque} T. Clausen et K. Tent aboutissent à un résultat similaire (voir Proposition 6.2 de \cite{CT}).\end{Remarque}
Nous rappelons maintenant quelques faits concernant les groupes de type impair dont le rang de Prüfer est égal à un. 
\begin{fait} \label{RP=1}\cite[Proposition 1]{DJ1} Soit $G$ un groupe rangé connexe de type impair, de rang de Prüfer égal à un. Alors il y a trois possibilités pour le type d'isomorphisme des 2-sous-groupes de Sylow ($S^{\circ}$ est un 2-tore de rang 1):
\begin{enumerate}
\item $S=S^{\circ}\rtimes \langle i \rangle$ pour $i$ une involution agissant sur $S^{\circ}$ par inversion (type $\mbox{PGL}_2 (\mathbb{C})$).
\item $S=S^{\circ}\cdot \ \langle w \rangle$ pour $w$ un élément d'ordre 4 agissant sur $S^{\circ}$ par inversion. On a de plus que $w^2$ est l'involution de $S^{\circ}$ (type $\mbox{SL}_2 (\mathbb{C})$).
\item $S=S^{\circ}$.
\end{enumerate}
\end{fait}
\begin{fait}
\label{Steiberg} \cite{Del2}
Soit $G$ un groupe rangé, $U_p^{\perp}$, et soit $x\in G$ un $p$-élément tel que $x^{p^n}\in Z(G)$. Alors l'exposant du groupe $C_G(x)/C_G(x)^{\circ}$ divise $p^n$.
\end{fait}
La proposition suivante caractérise complètement la structure de la 2-torsion des quasi groupes de Frobenius rangés connexes de degré pair et de type impair:
\begin{fait}\label{2-torsion quasi-frobenius} \cite[Proposition 1]{DW} Soit $C<G$ un quasi groupe de Frobenius rangé connexe tel que $C$ est définissable. Supposons que $G$ est de degré pair et de type impair. Alors $G$ a les mêmes 2-sous-groupes de Sylow que $\mbox{PGL}_2(\mathbb{C})$. De plus, on a $N_G(C)=C\rtimes \langle k \rangle$ avec $k$ une involution agissant par inversion. En particulier, $C$ est un groupe abélien. \end{fait} 
On remarque que dans ce cas $N_G(C)\setminus C$ est composé d'involutions agissant par inversion. En effet, $ckck=cc^{-1}=1$; soit $j$ une telle involution, on a $C^j=C$ et l'action est par inversion car $C_j\neq C$. 

On peut déterminer de la même manière la structure des 2-sous-groupes de Sylow pour les groupes de Frobenius rangés connexes de type impair. En effet, le rang de Prüfer est égal à un (Fait \ref{involution Frobenius}).
\begin{proposition}\label{Sylow Frobenius}
Soit $C<G$ un groupe de Frobenius rangé connexe de type impair. Alors les 2-sous-groupes de Sylow de $G$ sont soit connexes, soit du type $\mbox{SL}_2(\mathbb{C})$
\end{proposition}
\begin{proof}
Soit $S$ un 2-sous-groupe de Sylow et $S^{\circ}$ sa composante connexe; on peut l'étendre à un tore décent maximal $T$. Il existe donc un conjugué $C^g$ tel que $T\leq C^g$. On note $i$ l'involution contenue dans $C^g$.
Soit $j\neq i$ une involution et soit $C_j$ le conjugué de $C$ qui la contient. Si $j$ agissait par inversion sur $S^{\circ}$ alors $j$ commuterait avec $i$. Par conséquent, $1\neq i\in C_j$, contradiction. Donc les 2-sous-groupes de Sylow ne sont pas du type $\mbox{PGL}_2(\mathbb{C})$, et on conclut par le Fait \ref{RP=1}. \end{proof}
Pour les quasi groupes de Frobenius rangés connexes de degré impair $> 1$  et de type impair, nous n'avons pas de contrôle sur le rang de Prüfer. Cependant, les involutions forment une classe de conjugaison.
\begin{fait} \cite[Proposition 1]{DW}\label{conjugaison involutions}
Soit $C<G$ un quasi groupe de Frobenius rangé connexe tel que $C$ est définissable. Supposons que $G$ est de type impair et de degré impair. Alors $N_G(C)$ est un sous-groupe fortement inclus. En particulier, les involutions forment une classe de conjugaison dans $G$.
\end{fait}
\begin{Remarque}\label{conjugaison involutions bis}
Soit $C<G$ un quasi groupe de Frobenius rangé connexe tel que $C$ est définissable. Supposons que $G$ est de degré pair et de type impair. Alors les involutions sont conjuguées. En effet, le rang de Prüfer est égal à un.
\end{Remarque}
\subsection{Le cas résoluble}
Avant d'étudier le cas résoluble, rappelons un certain nombre de résultats utiles sur les sous-groupes de Carter, \textit{i.e.}, les sous-groupes définissables connexes nilpotents presque-autonormalisants.
\begin{fait}\label{Carter}
\begin{enumerate}
\item \cite[Théorème 3.11]{FJ} Soit $G$ un groupe rangé. Alors $G$ contient un sous-groupe de Carter.
\item \cite[Théorème 3.1]{Jal1} Soit $G$ un groupe rangé connexe. Alors les sous-groupes de Carter généreux sont conjugués.
\item \cite[Théorème 29]{Wag2} Soit $G$ un groupe rangé connexe résoluble. Alors ses sous-groupes de Carter sont conjugués.
\item \cite[Corollaire 2.10]{Fre5} Soit $G$ un groupe rangé. Alors tout tore décent de $G$ est contenu dans un sous-groupe de Carter.
\item \cite[Théorèmes 1.1 et 1.2]{Fre3} Soit $G$ un groupe rangé connexe résoluble. Alors tout sous-groupe $H$ contenant un sous-groupe de Carter est définissable, connexe et autonormalisant.
\item \cite[Théorème 3.11]{FJ} Soit $G$ un groupe rangé résoluble. Alors ses sous-groupes de Carter sont généreux et génériquement disjoints de leurs conjugués.
\item\cite[I.Corollaire 8.30]{ABC08} Soit $G$ un groupe rangé connexe résoluble. Alors $G=G'Q$ où $Q$ est un sous-groupe de Carter.
\end{enumerate}
\end{fait}
\begin{Lemme} \label{quasi-frobenius résoluble} (à comparer avec \cite[Lemme 11.32]{BN})
Soit $C < G$ un quasi groupe de Frobenius rangé connexe tel que $C$ est définissable. Supposons que $G$ est résoluble. Alors $G$ n'est pas nilpotent et $G=G'\rtimes C$ avec $C$ abélien. De plus, on peut interpréter un corps $K$ tel que $C\hookrightarrow K^{\times}$.
\end{Lemme}
\begin{proof}
Soit $H$ un sous-groupe définissable connexe nilpotent normal. Alors $H\cap C = \{ 1 \}$. En effet, soit $1\neq c\in H\cap C$; alors $Z(H)^{\circ}$ (infini car $H$ est infini connexe nilpotent) commute avec $c$ et donc $Z(H)^{\circ}\leq C$. Mais alors $H\leq C_G(Z(H)^{\circ})^{\circ}\leq C$ et donc $C=G$, contradiction. Le groupe $G$ est en particulier non nilpotent.
Le groupe $G'$ est un sous-groupe définissable connexe nilpotent normal et donc $G'\cap C=\{ 1 \} $. Par conséquent, $C$ est abélien et c'est un sous-groupe de Carter (il est en particulier autonormalisant, Fait \ref{Carter}). D'après le Fait \ref{Carter}, on a donc $G=G'\rtimes C$. 

Soit $B\leq G'$ un sous-groupe définissable connexe abélien $C$-minimal. On a bien $C_C(B)=\{1 \}$ : soit $1\neq c\in C_C(B)$; alors $B\leq C_G(c)^{\circ}\leq C$, contradiction. On peut donc appliquer \cite[Théorème 9.1]{BN} pour conclure.
\end{proof}
\begin{corollaire}\label{quasi-Frobenius résoluble}
Soit $C<G$ un quasi groupe de Frobenius rangé connexe tel que $C$ est définissable. Supposons que $G$ est résoluble. Alors $G$ est un groupe de Frobenius.
\end{corollaire}
\begin{Lemme} \label{noyau frobenius}
Soit $G=U\rtimes C$ un groupe rangé où $U$ est définissable connexe nilpotent et chaque $c\neq 1$ agit sur $U$ en fixant au plus un nombre fini de points. Alors $U\setminus \{1\}=G\setminus \bigcup_G C^g$.
\end{Lemme}
\begin{proof}
Par l'exercice 10.a p.98 de \cite{BN}, on a : $U=[U,c]$ pour $c\in C$. Soit $g=uc\in G\setminus U$ avec $u\in U$ et $c\in C$. Soit $u_1\in U$ tel que $u=[u_1^{-1},c^{-1}]$; alors $g^{u_1}\in C$.\end{proof}

Récapitulons maintenant les résultats obtenus dans le cas résoluble.
\begin{theoreme}\label{recap résoluble}
Soit $C<G$ est un quasi groupe de Frobenius rangé connexe tel que $C$ est définissable. Supposons que $G$ est résoluble, alors :
\begin{enumerate}
    \item $C<G$ est un groupe de Frobenius.
    \item $C$ est abélien.
    \item $G$ n'est pas nilpotent.
    \item $G'=G\setminus \{\bigcup_G C^g\}\cup \{1\}$.
    \item $G=G'\rtimes C$. 
    \item Il existe un corps définissable $K$ tel que $C\hookrightarrow K^{\times}$ et $RM(C)\leq RM(G')$.
\end{enumerate}
\end{theoreme}
\subsection{Le groupe de Weyl}
Cette section est une digression concernant la notion de groupe de Weyl au sein des quasi groupes de Frobenius rangés connexes. Une façon d'analyser cette classe de groupes consiste à s'intéresser à $N_G(C)/C$, une forme de \textit{groupe de Weyl} . En effet, pour un groupe algébrique réductif, le groupe de Weyl désigne $N_G(T)/C_G(T)$ où $T$ est un tore (algébrique) maximal. Or, pour les  groupes de Frobenius algébriques connexes, le complément de Frobenius est précisément un tore maximal (on renvoie à \cite[Proposition 1]{Her}). De même dans $\mbox{PGL}_2(K)$ qui est l'exemple paradigmatique de quasi groupe de Frobenius algébrique de degré pair, le quasi complément de Frobenius est un tore (algébrique) maximal. 

En général dans le contexte rangé, plusieurs notions de groupes de Weyl peuvent être envisagées. Si on suppose que le quasi complément de Frobenius est résoluble, alors elles coïncident toutes pour les quasi groupes de Frobenius rangés connexes de type impair. Nos résultats s'inspirent largement de \cite{ABFa} qui étudie les groupes de Weyl des groupes simples minimaux connexes rangés, \textit{i.e.}, des groupes simples rangés tels que les sous-groupes définissables connexes propres sont résolubles.

Il nous faut d'abord décrire la structure des sous-groupes de Carter dans les quasi groupes de Frobenius rangés connexes. 

\begin{Lemme}\label{Carter quasi-Frobenius}
Soit $C<G$ un quasi groupe de Frobenius rangé connexe tel que $C$ est définissable et résoluble. Alors $C$ contient un sous-groupe de Carter de $G$ généreux dans $G$. De plus, un sous-groupe de Carter de $G$ est soit généreux contenu dans un conjugué de $C$, soit il n'est pas généreux, intersecte trivialement $C$ et ses conjugués, et il ne contient pas de torsion divisible.
\end{Lemme}
\begin{proof}
D'après le Fait \ref{Carter}, le groupe $C$ contient un sous-groupe de Carter $Q$. Or, d'après le Fait \ref{Carter}, les sous-groupes de Carter d'un groupe connexe résoluble sont autonormalisants, conjugués et généreux. Puisque $C$ est généreux et connexe, par le Fait \ref{générosité}, le groupe $Q$ est un sous-groupe de Carter de $G$ généreux dans $G$ car $N_G(Q)^{\circ}=N_C(Q)^{\circ}=Q$. En particulier, tous les sous-groupes de Carter de $C$ sont des sous-groupes de Carter généreux de $G$.

Si $Q_1$ est un sous-groupe de Carter généreux de $G$, alors $Q_1\cap C$ (quitte à conjuguer $C$) est non-trivial. Puisque $Z(Q_1)^{\circ}\neq \{1\}$, on a $Q_1\leq C$. Inversement, si $Q_1$ est un sous-groupe de Carter qui intersecte $C$ non-trivialement, alors $Q_1\leq C$ (condition du normalisateur) et c'est un sous-groupe de Carter généreux.
Si $Q_1$ intersecte trivialement $C$ et ses conjugués, il ne contient pas de torsion divisible.\end{proof}
\begin{proposition} \label{groupe de Weyl quasi-frobenius} (à comparer avec \cite[Proposition 3.2]{ABFa})
Soit $C<G$ un quasi groupe de Frobenius rangé connexe tel que $C$ est définissable. Supposons que $G$ est de type impair. Alors on a $N_G(C)=CN_G(T)$ où $T$ est un tore décent maximal de $C$.
Si de plus $C$ est résoluble alors $W(G)=N_G(T)/C_G(T)\simeq N_G(Q)/Q\simeq N_G(C)/C$ où $Q$ est un sous-groupe de Carter généreux de $G$.
\end{proposition}
\begin{proof}
Soit $T$ un tore décent maximal de $C$. D'après le Fait \ref{tore décent}, les tores décents maximaux de $C$ sont conjugués dans $C$. On applique un argument de Frattini au groupe $N_G(C)$ : soit $n\in N_G(C)$ on a $T^{n}\leq C$ et par conjugaison des tores décents maximaux, il existe $c\in C$ tel que $c^{-1}n\in N_G(T)$, par conséquent $N_G(C)=CN_{N_G(C)}(T)$. Mais $N_G(T)\leq N_G(C)$, par conséquent, $N_G(C)=CN_G(T)$.

Supposons maintenant que $C$ est résoluble. Soit $Q$ un sous-groupe de Carter généreux de $G$. Quitte à conjuguer, on peut supposer que $T\leq Q\leq C$. On rappelle que les sous-groupes de Carter de $C$ sont conjugués et autonormalisants. Un argument de Frattini appliqué à $N_G(C)$ nous donne : $N_G(C)=CN_G(Q)$ (car $N_{N_G(C)}(Q)=N_G(Q)$). Par conséquent, $N_G(C)/C\simeq N_G(Q)/(C\cap N_G(Q))=N_G(Q)/Q$. De plus, $Q\leq C_G(T)\leq C$; donc $C_G(T)$ est définissable connexe résoluble. Il suit que $Q$ est autonormalisant dans $C_G(T)$. Cette fois-ci un argument de Frattini appliqué à $N_G(T)$ nous donne $N_G(T)=C_G(T)N_{G}(Q)$ (car $N_G(Q)\leq N_G(T)$ puisque $T$ est caractéristique dans $Q$).  Puisque $Q$ est autonormalisant dans $C_G(T)$, on a $N_G(T)/C_G(T)\simeq N_G(Q)/Q$. Finalement, $W(G)\simeq N_G(Q)/Q\simeq N_G(C)/C$.
\end{proof}
On obtient de façon similaire la proposition suivante (pour plus détails, nous renvoyons au chapitre 4 de  la thèse doctorale de l'auteur) :
\begin{proposition} \label{groupe de weyl p-tore} (à comparer avec \cite[Corollaire 3.4]{ABFa})
Soit $C<G$ un quasi groupe de Frobenius rangé connexe tel que $C$ est définissable et résoluble. Alors pour un $p$-tore maximal non-trivial $S$, on a $N_G(S)/C_G(S)\simeq W(G)$.
\end{proposition}
\begin{corollaire} \label{RP QFI}
Soit $C<G$ un quasi groupe de Frobenius rangé connexe tel que $C$ est définissable et résoluble. Supposons que $G$ est de degré impair > 1 et de type impair. Alors le $2$-rang de Prüfer est supérieur ou égal à 2. 
\end{corollaire}
\begin{proof}
Supposons que le sous-groupe $C$ contient un 2-tore maximal $T_2$ de rang de Prüfer égal à un. D'après la proposition \ref{groupe de weyl p-tore}, $N_G(T_2)/C_G(T_2)\simeq N_G(C)/C\neq \{ 1 \}$. Or, d'après le Fait 29 de \cite{DJ1}, $Aut(T_2)\simeq Z_{2}^{\times}$ et le seul automorphisme d'ordre fini non-trivial est l'inversion. Le groupe $N_G(C)/C$ serait d'ordre 2, contradiction. \end{proof}
\subsection{\` A la recherche d'un sous-groupe de Borel fortement standard}
L'étude des involutions et des translations constitue un pan essentiel de l'analyse des groupes  de type impair. Rappelons que A. Deloro et J. Wiscons montrent un résultat très important sur la situation des translations dans notre contexte (Théorème \ref{transpositions QFP}) :

Soit $C<G$ un quasi groupe de Frobenius rangé connexe tel que $C$ est définissable et connexe. Supposons que $G$ est de type impair. Alors $\bigcup_{g\in G} C^g$ ne contient pas toutes les translations.
\begin{Remarque}
Si $C<G$ est un groupe de Frobenius connexe de type impair, alors aucune translation n'est contenue dans $\bigcup_{g\in G} C^g$ : soit $x=ij\in C^g$; alors $i,j\in N_G(C^g)=C^g$, et donc $i=j$ (car le complément de Frobenius contient une unique involution).
\end{Remarque}
Sur la base de ce théorème, nous allons montrer l'existence d'un sous-groupe de Borel non-nilpotent qui intersecte $C$ (quitte à conjuguer) dans une partie infinie mais qui n'est pas contenu dans $C$. Idéalement, on conjecture qu'un tel sous-groupe de Borel devrait contenir tout $C$.
\begin{definition} Soit $C<G$ un quasi groupe de Frobenius rangé connexe tel que $C$ est définissable. Soit $B$ un sous-groupe de Borel. Alors :
\begin{itemize}
    \item $B$ est faiblement standard si $1<(B\cap C)< B$.
    \item $B$ est standard si $1<(B\cap C)< B$ et $B\cap C$ contient des involutions.
    \item $B$ est fortement standard si $C < B$.
\end{itemize}
 \end{definition}
Autrement dit, nous allons montrer que $G$ contient un sous-groupe de Borel faiblement standard. D'après le Théorème \ref{recap résoluble}, un sous-groupe de Borel $B$ faiblement standard n'est pas nilpotent, c'est un groupe de Frobenius tel que $B=B'\rtimes (B\cap C)$ avec $B\cap C\hookrightarrow K^{\times}$ pour un corps $K$ définissable. Si $B$ est standard, le sous-groupe $(B\cap C)$ contient des involutions qui inversent donc le groupe abélien $B'$. Finalement, un sous-groupe de Borel fortement standard $B$ s'écrit sous la forme $B=B'\rtimes C$.
Remarquons que dans certaines configurations de type CiBo le quasi complément $C$ est précisément un sous-groupe de Borel. L'existence d'un sous-groupe de Borel fortement standard est donc une hypothèse très forte qui permettrait d'éliminer des configurations pathologiques.

On utilise le terme `` standard '' car dans le cas d'un quasi groupe de Frobenius rangé connexe de degré pair, nous cherchons à identifier un sous-groupe de Borel dont la structure est similaire à celle des sous-groupes de Borel de $\mbox{PGL}_2(K)$ à savoir isomorphes au groupe strictement 2-transitif $K^{+}\rtimes K^{\times}$. Remarquons que dans $\mbox{PGL}_2(\mathbb{C})$ un sous-groupe de Borel est généreux mais il n'est pas génériquement disjoint de ses conjugués. En effet, soient $B\neq B^g$ deux sous-groupes de Borel; ils s'intersectent dans un tore maximal $T$ qui est autocentralisant dans $G$ (groupe réductif). Or, $C_B(T)=T$ est généreux dans $B$, par conséquent $RM(B\setminus \bigcup_{G/N_G(B)} B^g)<RM(B)$.
Pour les groupes de Frobenius connexes, nous cherchons plutôt à montrer qu'ils sont résolubles.

Tout d'abord, voici un lemme inspiré de la démonstration du Théorème B de \cite{DW}.
\begin{Lemme}\label{transposition complément} Soit $C<G$ un quasi groupe de Frobenius rangé connexe tel que $C$ est définissable. Supposons que $G$ est de type impair. Soient $i, j$ deux involutions. Si $C$ est $\langle i, j \rangle$-invariant, alors $ij\in \bigcup_G C^g$. \end{Lemme}
\begin{proof}
Si la configuration est de degré impair, alors $i, j\in C$. Reste le cas de degré pair. Alors $i, j\in N_G(C)$. Si les deux involutions sont dans $C$, $1=ij\in C$ (Fait \ref{2-torsion quasi-frobenius}). Si une seule est dans $C$, alors elle est centralisée par l'autre : $ij$ est une involution donc elle appartient à $\bigcup_G C^g$. Enfin, si $i, j\in N_G(C)\setminus C$, alors $ij\in C_G(C)$ et donc ou bien $ij\in C$, ou bien $ij$ est une involution contenue dans un conjugué de $C$.
\end{proof}
Nous pouvons maintenant énoncer le théorème :
\begin{theoreme} \label{borel standard''}
Soit $C<G$ un quasi groupe de Frobenius rangé connexe tel que $C$ est définissable. Supposons que $G$ est de type impair. Alors il existe un sous-groupe $A$ définissable connexe abélien inversé par (au moins) deux involutions $i\neq j$, tel que $A\cap \bigcup_G C^g=\{ 1 \}$, et maximal pour ces propriétés. De plus, on a les propriétés suivantes :
\begin{enumerate}
\item  $C_G(A)^{\circ}=A$.
\item $A$ est TI.
\item  Soit $N_G(A)=N$. Alors $A$ n'est pas presque-autonormalisant et $N^{\circ}=A\rtimes (N^{\circ}\cap C)$. 
\end{enumerate}
\end{theoreme}
\begin{proof}
Par le Théorème \ref{transpositions QFP}, on trouve une translation $x=ij$ qui n'est pas contenue dans $\bigcup_G C^g$, pour deux involutions $i\neq j$. On considère $A=C_G(x)^{\circ}$ (qui est infini par le Fait \ref{centralisateur infini}). Remarquons qu'il s'agit d'un sous-groupe $\langle i,j \rangle$-invariant tel que $C_A(i)^{\circ}=C_G(i,j)^{\circ}=C_A(j)^{\circ}$. Supposons que le groupe  $\langle i, j \rangle$-invariant $C_A (i)^{\circ}$ est infini. On a donc  $C_A (i)^{\circ}\leq C_G(i)^{\circ}\leq  C_i$. Puisque $C_i$ est TI, il est aussi $\langle i, j\rangle $-invariant et donc $x\in \bigcup_G C^g$ (Lemme \ref{transposition complément}), contradiction.
Par conséquent, le groupe $C_A (i)=C_A(j)$ est fini et donc $A$ est un groupe définissable connexe abélien inversé par $i$ et $j$.  De plus, $A\cap \bigcup_G C^g=\{ 1 \}$; le groupe $A$ est en particulier uniquement 2-divisible.

Il existe donc un sous-groupe définissable connexe abélien inversé par au moins deux involutions et d'intersection triviale avec les conjugués de $C$. On prend $A$ maximal pour ces propriétés et on pose $i\neq j$ deux involutions l'inversant. Montrons tout d'abord $C_G(A)^{\circ}=A$. On considère $C_G(A)^{\circ}$; c'est un groupe $\langle i,j \rangle$-invariant, inversé par $i,j$. De plus, $C_G(A)^{\circ}\cap \bigcup_G C^g=\{ 1 \}$. Par maximalité de $A$ (et puisque $A$ est abélien), on a $C_G(A)^{\circ}=A$. Montrons maintenant que le groupe $A$ est TI. Supposons par l'absurde qu'il existe un conjugué $A'\neq A$ et $1\neq x'\in A'\cap A$. L'élément $x'$ est une translation inversée par $i$. De plus, $A,A'\leq C_G(x')^{\circ}=A_1$. En raisonnant comme précédemment, on obtient que le groupe $A_1$ est abélien inversé par $i$ et $j$. Ainsi, $A_1=C_G(A)^{\circ}=A=A'$, contradiction. Donc $A$ est TI.

Le sous-groupe $A$ n'est pas presque-autonormalisant : si $N_G(A)^{\circ}=A$, alors le groupe $A$ est généreux et donc intersecte non-trivialement un conjugué de $C$, contradiction. On pose $N=N_G(A)$. Remarquons que $i,j\in N$. Si les groupes $C_N (i)$ et $C_N (j)$ étaient finis alors $i,j$ inverseraient le groupe abélien $N^{\circ}$. En particulier, $N^{\circ}=C_G(A)^{\circ}=A$, contradiction. Par conséquent, $(N^{\circ}\cap C_i)<N^{\circ}$ est un quasi groupe de Frobenius connexe avec complément définissable (Lemme \ref{conjugaison complément Frobenius}). On lui applique la Proposition \ref{scission quasi-frobenius} pour obtenir que $N^{\circ}=A\rtimes(C_i\cap N^{\circ})$. \end{proof}
Si $C$ est résoluble, pour obtenir un sous-groupe de Borel faiblement standard, il suffit de considérer un sous-groupe de Borel $B$ contenant le groupe connexe résoluble $N^{\circ}$. 
\begin{Remarque} On rappelle qu'un sous-groupe de Borel faiblement standard $B$ s'écrit sous la forme $B=B'\rtimes (B\cap C)$, avec $(B\cap C)\hookrightarrow K^{\times}$, pour un corps $K$. Par conséquent, si $B$ contient un $p$-groupe unipotent, alors sous l'hypothèse de l'infinité des $p$-nombres premiers de Mersenne, le Fait \ref{p-Mersennes} nous permet de conclure que $B$ est un sous-groupe de Borel standard.
\end{Remarque}
Comment obtenir un sous-groupe de Borel fortement standard ? Dans \cite{FJ}, les auteurs considèrent plusieurs conjectures de généricité, notamment la \textit{générosité des sous-groupes de Borel} (GB) et \textit{la générosité des sous-groupes de Carter} (GC) (dans le deuxième cas, par \cite{Jal1}, cela implique la conjugaison des sous-groupes de Carter). On peut envisager ces hypothèses de générosité comme la reformulation dans le contexte rangé d'hypothèses topologiques de densité.

L'hypothèse (GB) est suffisante pour démontrer l'existence d'un sous-groupe de Borel fortement standard si le quasi complément de Frobenius est nilpotent. Néanmoins, dans les configurations de type CiBo, il existe des sous-groupes de Borel qui ne sont pas généreux (voir \cite{CJ}, par exemple le Lemme 5.11). 

\begin{Lemme} \label{critère pour borel standard}
Soit $C<G$ un quasi groupe de Frobenius rangé connexe tel que $C$ est définissable et nilpotent. Soit $B$ un sous-groupe de Borel tel que $1<B\cap C<B$. Si $B$ est généreux, alors $C<B$.
\end{Lemme}
\begin{proof} 
 D'après le Lemme \ref{conjugaison complément Frobenius}, $(B\cap C)<B$ est un quasi groupe de Frobenius connexe avec complément définissable. Notamment, $B\cap C$ est généreux dans $B$, qui l'est dans $G$; comme en outre $B$ est connexe, par transitivité $B\cap C$ est généreux dans $G$. Notamment, $N_G(B\cap C)^{\circ}=(B\cap C)$. Par nilpotence de $C$ et croissance des normalisateurs, $C<B$.
\end{proof}
\subsection{Théorèmes de classification}
Nous pouvons maintenant démontrer nos deux théorèmes de classification :
\begin{enumerate}
    \item Théorème \ref{classification QFP} : 
    
    Soit $C<G$ un quasi groupe de Frobenius connexe rangé tel que $C$ est définissable. Supposons que $G$ est de degré pair et de type impair. Si les sous-groupes de Borel sont généreux, alors $G$ est isomorphe à $\mbox{PGL}_2(K)$ pour un corps $K$ (rangé) algébriquement clos de caractéristique différente de 2.
    \item Théorème \ref{classification frobenius} : 
   
    Soit $C<G$ un groupe de Frobenius connexe rangé de type impair avec $C$ résoluble et $G$ non résoluble. Si les sous-groupes de Borel sont généreux, alors $C$ n'est pas nilpotent.
\end{enumerate}
On démontrera également les théorèmes précédents en remplaçant l'hypothèse de générosité des sous-groupes de Borel par celle de l'existence d'un sous-groupe de Borel fortement standard.
Nous formulerons certains lemmes intermédiaires sous l'hypothèse de l'existence d'un sous-groupe de Borel faiblement standard (ou standard). Nous souhaitons dans la mesure du possible éviter d'utiliser toute la force de l'hypothèse de l'existence d'un sous-groupe de Borel fortement standard.

Pour les quasi groupes de Frobenius rangés connexes de degré pair, nous savons que le quasi complément de Frobenius est abélien. Néanmoins, rappelons qu'il pourrait s'agir d'un sous-groupe de Borel (une situation fortement non-algébrique). On retrouve là l'un des contre-exemples  ($\mbox{CiBo}_2$) minimaux potentiels à la Conjecture de Cherlin-Zilber pour les groupes de petit rang de Prüfer (voir \cite{DJ3}). Une résolution complète de la conjecture $A_1$ de \cite{DW} permettrait donc l'élimination de cette configuration pathologique.

L'identification de $\mbox{PGL}_2(K)$ se fait \textit{via} l'existence d'une $(B,N)$-paire scindée spéciale de rang 1. Dans tout groupe algébrique réductif, un sous-groupe de Borel $B$ et $N_G(T)$ où $T$ est un tore maximal contenu dans $B$ forment une $(B,N)$-paire. Le rang de la $(B,N)$-paire correspond au nombre d'involutions qui engendrent le groupe $N/B\cap N$.  Plus spécifiquement, dans $\mbox{PGL}_2(K)$, une telle $(B,N)$-paire est scindée de rang 1 (on a en effet $B=B_u\rtimes T$ où $B_u$ est le radical unipotent et $T$ le tore algébrique maximal). De plus, dans ce cas, un tore maximal $T$ est un groupe TI presque-autonormalisant tel que $N_G(T)/T$ est d'exposant 2.

Nous nous inspirons de la stratégie employée dans \cite{DJ3} (voir en particulier la Proposition 3). Nous renvoyons également à \cite{Del} et \cite{CJ}.

Une $(B,N)$-paire scindée spéciale de rang 1 est un quadruplet de groupes $B$, $N$, $H=N\cap B$ et $U$ un groupe nilpotent avec $H$ normal dans $N$, qui satisfait les axiomes suivants :
\begin{enumerate}
\item  $G=\langle B,N \rangle$.
\item $[N:H]=2$.
\item Pour tout $k\in N\setminus H$, on a $H=B\cap B^k$, $G=B\cup BkB$, et $B^k\neq B$.
\item $B=U\rtimes H$.
\end{enumerate}

Les deux théorèmes suivants permettent d'identifier des groupes algébriques \textit{via} des $(B,N)$-paires dans le contexte rangé.
\begin{fait} \label{BN paire caractéristique nulle} (adapté de \cite[Théorème 2.1]{DMT}) Dans un univers rangé, soit $(B,N,U,H)$ une $(B,N)$-paire définissable scindée spéciale de rang 1 telle que que $U$ est abélien divisible sans torsion. Alors $G\simeq \mbox{PGL}_2(K)$.\end{fait} 
\begin{fait} \label{BN paire caractéristique positive} (adapté de \cite[Théorème 1.2]{Wis}) Dans un univers rangé, soit $(B,N,U,H)$ une $(B,N)$-paire définissable scindée spéciale de rang 1. Supposons de plus que $U$ est un groupe abélien contenant un élément d'ordre $p>2$. Soient $L=\bigcap_{G}B^g$ et $M=\langle U^g : g\in G \rangle L$; si $H\neq L$ et si $(H\cap M)/L$ ne contient pas de $p$-groupe abélien élémentaire, alors $G/L\simeq \mbox{PGL}_2(K)$. \end{fait}
\subsubsection{Sous-groupe de Borel standard}
\begin{nota}
Soit $B$ un sous-groupe définissable et soit $k$ une involution qui ne normalise pas $B$. On pose :
\[ T_B(k)=\{b\in B : b^k=b^{-1} \}.\]
\end{nota}
\begin{Lemme}\label{involution générique}
Soit $C<G$ un quasi groupe de Frobenius rangé connexe tel que $C$ est définissable et résoluble. Supposons que $G$ n'est pas résoluble et qu'il est de type impair. Soit $B$ un sous-groupe de Borel standard  et soit $i$ une involution dans $B\cap C$ telle que $N_G(B)\cap i^G=B\cap i^G$. 
Alors $i^G\setminus N_G(B)$ et \[K_B=\{ k\in i^G\setminus N_G(B) : RM(T_B(k))\geq RM(B)-RM(C_G(i))\}\] sont génériques dans $i^G$. 
\end{Lemme}
\begin{proof}
Les involutions de $G$ forment une classe de conjugaison $i^G$ (Fait \ref{conjugaison involutions} et Remarque \ref{conjugaison involutions bis}). Par hypothèse, $N_G(B)\cap i^G=B\cap i^G$. Supposons que $i^G\setminus N_G(B)$ n'est pas générique dans $i^G$, \textit{i.e.}, $i^G \cap N_G(B)=i^G\cap B$ est large dans $i^G$. Le groupe $N=\bigcap_G B^g=\bigcap_{1\leq i \leq n}B^{g_i}$ est un groupe définissable normal contenant des, donc toutes les involutions (car $i^G\setminus \bigcup_G B^g=i^G\setminus \bigcup_{1\leq i \leq n}B^{g_i}$ est non-générique); c'est en particulier un sous-groupe infini.

Si $H=C\cap N^{\circ}$ est trivial, alors par la Proposition \ref{scission quasi-frobenius}, $G=A\rtimes C$ pour un sous-groupe définissable abélien $A$ et $G$ est résoluble, contradiction.
Le sous-groupe $H<N^{\circ}$ est donc un quasi groupe de Frobenius connexe avec complément définissable car $N^{\circ}\not\leq C$ (Lemme \ref{conjugaison complément Frobenius}). Mais pour tout $g\in G$, il existe $n\in N^{\circ}$ tel que $H^g=H^n$. Par conséquent, $G=N^{\circ}\cdot N_G(H)$. Mais $N_G(H)\leq N_G(C)$. Or, $N_G(C)$ est résoluble; le groupe $G$ est donc résoluble, contradiction.

Supposons maintenant que $i^G\setminus N_G(B)$ est générique dans $i^G$. Dans la suite, nous reprenons la démonstration de la Proposition 2 de \cite{DJ3} pour obtenir que \[K_B=\{ k\in i^G\setminus N_G(B) : RM(T_B(k))\geq RM(B)-RM(C_G(i)) \}\]  est aussi générique dans $i^G$.
On considère la fonction définissable $\phi : i^G\setminus N_G(B)\longrightarrow G/B$ qui envoie $k$ sur sa classe $kB$.
Le domaine est de rang $RM(i^G)=RM(G)-RM(C_G(i))$ et l'image de rang au plus $RM(G)-RM(B)$. La fibre $\phi^{-1}(\phi(k))$ pour $k$ générique dans $i^G$ est donc de rang au moins $RM(B)-RM(C_G(i))$. Mais si $k$ et $j$ appartiennent à la même fibre, on a $jB=kB$ et donc $kj\in T_B(k)$.

Par conséquent, pour $k$ une involution générique, on a \[RM(T_B(k))\geq RM(\phi^{-1} \phi(k))\geq RM(B)-RM(C_G(i)).\] \end{proof} 
\begin{Lemme}\label{Frobenius standard}
Soit $C<G$ un groupe de Frobenius rangé connexe de type impair avec $C$ résoluble et $G$ non résoluble. Supposons qu'il existe un sous-groupe de Borel standard $B$ tel que $B'$ est TI. Alors $T_B(k)$ est fini, pour $k\in K_B$; de même que $B\cap B^k$. De plus, $RM(B)\leq RM(C)$.
\end{Lemme}
\begin{proof} 
Soit $i\in C$ l'unique involution de $C$. Pour un sous-groupe $H$, on pose $I(H)=\{1\neq h\in H : h^2=1\}$.
On a $N_G(B)\cap i^G = B\cap i^G$. En effet, soit $j\in N_G(B)$; alors $C_B(j)$ n'est pas trivial car $B$ n'est pas nilpotent (Théorème \ref{recap résoluble}). Il suit que $1<(B\cap C_j)<B$ est un groupe de Frobenius connexe. D'après le Lemme \ref{conjugaison complément Frobenius},  $(B\cap C_j)$ et $(B\cap C_i)$ sont conjugués dans $B$ : le sous-groupe $B\cap C_j$ contient une involution qui est égale à $j$ et donc $j\in B$. Enfin, $I(N_G(B\cap C))\leq I(N_G(C))=I(C)=\{i\}$.
D'après le Lemme \ref{involution générique}, $i^G\setminus N_G(B)$ et $K_B=\{ k\in i^G\setminus N_G(B) : T_B(k)\geq RM(B)-RM(C_G(i)) \}$ sont génériques dans $i^G$. 

Soit $k\in K_B$; montrons que si $x\in B$ est inversé par $k$ alors $x$ est trivial. Supposons que $x\neq 1$.
Supposons d'abord que $x\in B\cap C$. Alors $k$ normalise $C_G(x)^{\circ}$ puis $B\cap C$, donc $k=i$, contradiction. De même pour $x\in (B\cap C)^b$ avec $b\in B$.
D'après le Théorème \ref{recap résoluble}, $x\in B'$ et donc $x^{-1}\in B'\cap(B')^k$. Puisque $B'$ est TI, on a $B'=(B')^k$. Comme $N_G(B')^{\circ}$ est résoluble (Proposition \ref{scission quasi-frobenius}), on obtient $B^k=B$, contradiction.  
Il suit que $B\cap B^k=\{1\}$ : en effet, si $x$ est dans l'intersection, alors $y=xx^{-k}\in B$ est inversé par $k$ donc $y=1$; ainsi $x=x^k$ et $x=1$ pour les mêmes raisons.

Si $RM(B)>RM(C)$, alors $RM(T_B(k))\geq RM(B)-RM(C)>0$, contradiction.
\end{proof}
\begin{corollaire}
Soit $C<G$ un groupe strictement 2-transitif rangé de type impair avec $C$ résoluble et $G$ non résoluble. Soit $B=N_G(C_G(ij))$ pour deux involutions $i\neq j$, avec $i\in C$. Alors $B$ est un sous-groupe de Borel et pour $k\in K_B(k)$, $T_B(k)$ et l'intersection $B^k\cap B$ sont finis.
\end{corollaire}
\begin{proof}
On remarque que $G$ est un groupe de Frobenius connexe.
Soient $i\neq j$ deux involutions avec $i\in C$.
Le groupe $N=N_G(C_G(ij))=C_G(ij)\rtimes N_C(C_G(ij)$ est un groupe strictement 2-transitif par la Proposition 11.51 de \cite{BN}. De plus, $C_G(ij)=iI\cap jI$ est un sous-groupe TI \cite[Lemme 11.50]{BN}. Mais par le Fait 7 et le Théorème 10 de \cite{ABW}, on a : $N=K^{+}\rtimes K^{\times}$, pour $K$ un corps interprétable algébriquement clos. 
Montrons que $N$ est un sous-groupe de Borel standard tel que $N'$ est TI. Tout d'abord, d'après le Théorème \ref{recap résoluble}, $N'=C_G(ij)$ et $N'=N\setminus \bigcup_N (C\cap N)^n\cup \{1\}$. Soit $B$ un sous-groupe de Borel contenant $N$. D'après le Théorème \ref{recap résoluble}, $B=B'\rtimes (B\cap C)$ et $B'=B\setminus \bigcup_B (B\cap C)^b\cup \{1\}$; les involutions $i$ et $j$ agissent par inversion sur le groupe abélien $B'\subseteq iI\cap jI=N'\subseteq B'$ et donc $B'= N'$. D'après la Proposition \ref{scission quasi-frobenius}, $N_G(B')^{\circ}$ est un quasi groupe de Frobenius connexe scindé et résoluble contenant $B$; il lui est donc égal par maximalité de $B$. Par conséquent, on a $B=N$. Il suffit d'appliquer le Lemme \ref{Frobenius standard}.
\end{proof}
 \subsubsection{Sous-groupe de Borel fortement standard}
\begin{Lemme} \label{B autonormalisant}
Soit $C<G$ un quasi groupe de Frobenius rangé connexe tel que $C$ est définissable et résoluble. Supposons que $G$ est de degré pair ou de degré un, et de type impair. Supposons qu'il existe un sous-groupe de Borel fortement standard $B$. Alors $B^k\neq B$ pour tout $k\in N_G(C)\setminus C$. De plus, $B$ est autonormalisant.

\end{Lemme}
\begin{proof} Soit $i$ l'involution de $C$; d'après le Théorème \ref{recap résoluble}, on a $B=B'\rtimes C_i$.
On rappelle que si $G$ est de degré pair alors $N_G(C_i)\setminus C_i$ consiste d'involutions (Fait \ref{2-torsion quasi-frobenius}); si $G$ est de degré un, $N_G(C_i)=C_i$.
Soit $k\in N_G(C_i)\setminus C_i$ une involution (le cas de degré un est évident). Supposons par l'absurde que $B^{k}=B$. Le groupe $C_B(k)^{\circ}$ est infini (sinon $B$ serait abélien et donc égal à $C$). D'après le Lemme \ref{conjugaison complément Frobenius}, $C_k<B$ et donc $k\in B$. Mais d'après le Théorème \ref{recap résoluble}, $B$ est un groupe de Frobenius et donc $k\in C_i$, contradiction.

Montrons maintenant que $B$ est autonormalisant. Soit $g\in N_G(B)\setminus B$; le groupe $C_i^g$ est un quasi complément de Frobenius de $B$ et par le Lemme \ref{conjugaison complément Frobenius}, il existe $b\in B$ tel que $C_i^{gb^{-1}} = C_i$. On a : $g\in B$ ssi $gb^{-1} \in B$. On peut donc supposer que $g\in N_G(C_i)$. Si $g\in C_i$ alors on conclut, sinon $g$ est une involution normalisant $C_i$ et $B$, le paragraphe précédent montre qu'alors $g=i\in C_i$.
\end{proof}
\begin{Lemme} \label{intersection torale}
Soit $C<G$ un quasi groupe de Frobenius rangé connexe tel que $C$ est définissable et résoluble. Supposons que $G$ n'est pas résoluble et qu'il est de type impair. On suppose qu'il existe un sous-groupe de Borel fortement standard $B$. Soit $i$ une involution de $B\cap C$ telle que $i^G\setminus B$ et $K_B=\{ k\in i^G\setminus B : RM(T_B(k))\geq RM(B)-RM(C_G(i))\}$ sont génériques dans $i^G$. Alors pour toute involution $k\in K_B$, le groupe $I_k=(B\cap B^k)^{\circ}$ est un conjugué de $C$ contenant une involution $j_k$.
\end{Lemme}
\begin{proof}
On a $B=B'\rtimes C$. Soit $k\in K_B$; on a tout d'abord $RM(T_B(k))\geq RM(B)-RM(C_G(i))=RM(B')$. Le sous-groupe $I_k$ est $k$-invariant. Supposons que $C_{I_k}(k)^{\circ}\leq C_k$ est infini : le groupe $B\cap C_k$ est un quasi complément de Frobenius de $B$ donc conjugué à $C$ et finalement $C_k< B$, contradiction.

Le sous-groupe $I_k\subseteq T_B(k)\subseteq B\cap B^k$ est donc abélien inversé par $k$ : il suit que $RM(I_k)=RM(T_B(k))\geq RM(B')$.
Si $I_k\subseteq (B\setminus B')\cup \{1\}=\bigcup_B C^b$ (Théorème \ref{recap résoluble}), alors il existe $1\neq x\in I_k\cap C^b$, pour $b\in B$. Donc par commutativité, $I_k\leq C_G(x)^{\circ}\leq C^b$. Mais $RM(C^b)=RM(C)\leq RM(B')$ (Théorème \ref{recap résoluble}), donc $C^b=I_k$. 
Sinon, $I_k\leq B'$. Puisque $k$ agit par inversion et $RM(I_k)\geq RM(T_k(B))\geq RM(B')$, on a $B'= (B^k)'$.
On considère $N_G(B')^{\circ}$ qui contient $B$ et $B^k$. Mais alors $N_G(B')^{\circ}$ est un quasi groupe de Frobenius connexe contenant un groupe définissable connexe normal $B'$ d'intersection triviale avec $C$ et donc par la Proposition \ref{scission quasi-frobenius}, $N_G(B')^{\circ}=A\rtimes C$ où $A$ est abélien. Le sous-groupe $N_G(B')^{\circ}$ est définissable connexe résoluble et contient $B$ et $B^k$, contradiction (car $k$ ne normalise pas $B$).
\end{proof}
\begin{corollaire} \label{scission standard Frobenius}
Soit $C<G$ un groupe de Frobenius connexe rangé de type impair avec $G$ non résoluble et $C$ résoluble. Alors $C$ est un sous-groupe de Borel.
\end{corollaire}
\begin{proof}
Supposons que $C$ est contenu strictement dans un sous-groupe de Borel $B$. Soit $i$ l'unique involution de $C$. D'après le Lemme \ref{involution générique}, $i^G\setminus N_G(B)$ et $K_B=\{ k\in i^G\setminus N_G(B) : RM(T_B(k))\geq RM(B)-RM(C_G(i))\}$ sont génériques dans $i^G$.
Mais d'après le Lemme \ref{intersection torale}, pour $k\in K_B$, le groupe $(B\cap B^k)^{\circ}=I_k$ est un conjugué de $C$. Mais on a alors $k\in N_G(I_k)$ et donc $k\in I_k\leq B$, contradiction car $k$ ne normalise pas $B$.
\end{proof}
\begin{corollaire}\label{scission standard (GB) Frobenius}
Soit $C<G$ un groupe de Frobenius connexe rangé de type impair avec $G$ non résoluble et $C$ résoluble. Si les sous-groupes de Borel sont généreux, alors $C$ n'est pas nilpotent.
\end{corollaire}
\begin{proof} 
D'après le Théorème \ref{borel standard''}, il existe un sous-groupe de Borel faiblement standard qui est généreux par hypothèse.
Si $C$ est nilpotent, alors par le Lemme \ref{critère pour borel standard}, on a $C<B$, contradiction.
\end{proof}
\begin{Remarque} Dans \cite{CT}, T. Clausen et K. Tent aboutissent  également à un résultat de classification des groupes de Frobenius connexes de type impair avec complément nilpotent. Ils font l'hypothèse de l'absence de mauvais corps en caractéristique nulle. 
\end{Remarque}
Nous pouvons désormais compléter la suite des lemmes précédents pour identifier $\mbox{PGL}_2(K)$ parmi les quasi groupes de Frobenius rangés connexes de degré pair.
\begin{Lemme} \label{décomposition de Bruhat}
Soit $C<G$ un quasi groupe de Frobenius rangé connexe tel que $C$ est définissable. Supposons que $G$ est de degré pair et de type impair. Supposons qu'il existe un sous-groupe de Borel fortement standard $B$.
Soit $i$ l'involution de $B\cap C$; supposons  de plus que $i^G\setminus B$ et $K_B=\{ k\in i^G\setminus B : RM(T_B(k))\geq RM(B)-RM(C_G(i))\}$ sont génériques dans $i^G$ et que d'autre part pour tout involution $k\in K_B$, $I_k=(B\cap B^k)^{\circ}$ est un conjugué de $C$ contenant une involution $j_k$. 

Alors pour tout $g\in G\setminus B$, on a $G=B\cup BgB$ et donc $G$ est un groupe 2-transitif (pour l'action sur $G/B$).
\end{Lemme}
\begin{proof}
On rappelle que par le Lemme \ref{B autonormalisant}, le groupe $B$ est autonormalisant.
Puisque le sous-groupe $I_k$ est inversé par $k$, on a de plus $RM(T_B(k))=RM(I_k)\geq RM(B)-RM(C)=RM(B')$. Mais $I_k<B$ est un quasi groupe de Frobenius connexe avec complément définissable et par conséquent on peut interpréter un corps $K$ tel que $RM(B')\geq RM(K^{+})=RM(K^{\times})\geq RM(I_k)\geq RM(B')$ et donc $RM(B')=RM(I_k)=RM(C)$. Cette égalité est vraie pour les conjugués de $B$ et les conjugués de $C$ qui y sont contenus. En particulier, $B'$ est soit un groupe abélien divisible sans torsion soit un $p$-groupe abélien élémentaire. 

Considérons la fonction définissable de $K_B$ vers $i^B$ qui à chaque $k\in K_B$ associe l'unique involution $j_k$ de $I_k$. De plus, si $j_k=j_\ell$ alors $B\cap B^{\ell}=B\cap B^k$ et donc $\ell$ normalise le conjugué $C_{j_k}=I_k$, \textit{i.e.}, $\ell\in C_G(j_k)$. Les fibres sont de rang $RM(C_G(i))$. On a donc :
\[ RM(G)-RM(C_G(i))=RM(i^G)=RM(K_B)\leq RM(i^B)+RM(C_G(i))\]
\[ =RM(i^B)+RM(C_B(i))=RM(B). \]

Finalement, $RM(G)\leq RM(B)+RM(C_i)=RM(B)+RM(B')$.

Soit $g\in G\setminus B$; montrons que $B'\cap B^g=\{1\}$. Supposons par l'absurde qu'il existe $1\neq x\in B'\cap B^g$. On peut distinguer deux cas : $x\in (B')^g$ ou bien $x\in C_k\leq B^g$. Si $x\in (B')^g\cap B'$, alors $x$ est une translation telle que $B', (B')^g\leq C_G(x)^{\circ}$. Le sous-groupe $C_G(x)^{\circ}$ est abélien et donc $B', (B')^g\leq N_G(B')^{\circ}$, mais $N_G(B')^{\circ}=B$ (Proposition \ref{scission quasi-frobenius}). On en déduit que $B'=B'^g$, puis $B=B^g$, contradiction. Si $x\in B'\cap C_k$, alors $C_k\leq B$ et $C_k=C^b$ pour $b\in B$. Puisque $B'=B\setminus \bigcup_B C^b\cup \{1\}$, on obtient une contradiction.

Par conséquent, les fibres de l'application de $B'\times B$ vers $B'gB$ qui à $(f,b)$ associe $fgb$ sont finies. On a donc $RM(B'gB)=RM(G)$. La double classe $BgB$ est générique.

Par conséquent, pour tous $g,g'\in G\setminus B$,  les doubles classes $BgB$ et $Bg'B$ sont génériques dans $G$ et donc par connexité de $G$, $g\in Bg'B$. On a donc $G=B\cup Bg'B$. Il suit que le groupe $G$ est $2$-transitif pour l'action de $G$ sur $G/B$. \end{proof}
\begin{Remarque}
Le calcul de rang mis en \oe uvre dans la démonstration précédente utilise pleinement l'hypothèse de l'existence d'un sous-groupe de Borel \textit{fortement} standard.
\end{Remarque}
On peut désormais conclure la démonstration du Théorème \ref{classification QFP} :

Soit $C<G$ un quasi groupe de Frobenius rangé connexe tel que $C$ est définissable. Supposons que $G$ est de degré pair et de type impair. Si les sous-groupes de Borel sont généreux, alors $G\simeq \mbox{PGL}_2(K)$, pour un corps $K$ algébriquement clos de caractéristique différente de deux.
\begin{proof}
D'après le Fait \ref{2-torsion quasi-frobenius}, le groupe $C$ est abélien et contient une unique involution $i$. De plus, le groupe $G$ n'est pas résoluble (Corollaire \ref{quasi-Frobenius résoluble}). Il existe un sous-groupe de Borel $B$ qui contient strictement $C_i$, \textit{i.e.}, un sous-groupe de Borel fortement standard. En effet, d'après le Théorème \ref{borel standard''}, il existe un sous-groupe de Borel $B$ faiblement standard et puisque $B$ est généreux par hypothèse, le Lemme \ref{critère pour borel standard} montre que $B$ est en fait un sous-groupe de Borel fortement standard. 

Nous allons montrer que $(B,N_G(C_i))$ induit une $(B,N)$-paire (définissable) scindée spéciale de rang 1.
D'après le Lemme \ref{B autonormalisant}, pour tout $k\in N_G(C_i)\setminus C_i$, $B^k\neq B$. De plus, $B$ est autonormalisant.

On rappelle qu'on pose $T_B(k)=\{b\in B : b^{-1}=b^k \}$ pour une involution $k$ qui ne normalise pas $B$.  D'après le Lemme \ref{involution générique}, l'ensemble $K_B=\{ k\in N_G(B)\setminus B : RM(T_B(k))\geq RM(B)-RM(C_G(i))\}$ est générique dans $i^G$. 

On a $B\cap N_G(C_i)=C_i$ car $C_i<B$ est un groupe de Frobenius connexe. On a donc $[N_G(C_i): (N_G(C_i)\cap B)]=2$ (Proposition \ref{2-torsion quasi-frobenius}). De plus, pour $k\in N_G(C_i)\setminus C_i$, on a $C_i\leq B\cap B^k$ et $(B\cap B^k)=C_i$ car $(B\cap B^k)^{\circ}$ est abélien inversé par $k$ et $N_B(C_i)=C_i$. Par ailleurs, $B=B'\rtimes C_i$ où $B'$ est abélien. On a donc vérifié l'axiome 4, l'axiome 2 et une partie de l'axiome 3.

Soit $k\in K_B$; alors d'après le Lemme \ref{intersection torale}, $I_k=(B\cap B^k)^{\circ}$ est un conjugué de $C_i$ contenant une involution $j_k$.

Finalement par le Lemme \ref{décomposition de Bruhat}, on obtient :

$G=B\cup BkB$ pour tout $k\in N_G(C_i)\setminus C_i$ et $G=\langle B,N \rangle$.

On a donc montré que $(B, N_G(C_i), C_i, B')$ forme une BN-paire scindée spéciale de rang $1$.

Or, $\bigcap_{g\in G} B^g$ est trivial. En effet, on a $B\cap B^k=C_i$ pour $k\in N_G(C_i)\setminus C_i$ et donc $\bigcap_{g\in G} B^g\leq C_i$. Donc si $\bigcap_{g\in G} B^g\neq \{1\}$, on trouve $G=N_G(\bigcap_{g\in G} B^g)^{\circ}\leq N_G(C_i)^{\circ}=C_i$, contradiction. 

Si $B'$ est un $p$-groupe abélien élémentaire, alors en remarquant que $U_p(C_i)=\{ 1 \}$ (on peut interpréter un corps $K$ tel que $C_i\leq K^{\times}$) on applique le Fait \ref{BN paire caractéristique positive} et on a $G/\bigcap_{g\in G} B^g\simeq \mbox{PGL}_2(K)$ et donc $G\simeq \mbox{PGL}_2(K)$.
Sinon, $B'$ est abélien divisible sans torsion et on applique le Fait \ref{BN paire caractéristique nulle}. Dans les deux cas, on obtient $G\simeq \mbox{PGL}_2(K)$. \end{proof}

En guise de conclusion, notons que tout l'enjeu est de montrer que le sous-groupe de Borel faiblement standard obtenu en §4.5 est un sous-groupe de Borel fortement standard. Comme en témoigne l'étude des configurations de type CiBo, il s'agit d'un problème difficile. Néanmoins, nous espérons que les résultats présentés dans cet article seront susceptibles de rendre possible une approche inductive de ce problème, au moins en caractéristique positive (où l'existence de mauvais corps est improbable).

\end{document}